\documentclass[12pt,a4paper,draft]{amsart}

\usepackage{latexsym}
\usepackage{ifthen}
\usepackage[leqno]{amsmath}
\usepackage{enumerate}
\usepackage{calc}
\usepackage{amstext,amsbsy,amsopn,amsthm,amsgen,amsfonts,amscd,amsxtra,upref}
\usepackage{mathrsfs}\usepackage{euscript}\usepackage{amssymb}

\swapnumbers
\theoremstyle{plain}
\newtheorem{Thm}{Theorem}[section]
\newtheorem{Lem}[Thm]{Lemma}
\newtheorem{Cor}[Thm]{Corollary}
\newtheorem{Pro}[Thm]{Proposition}
\newtheorem{Prp}[Thm]{Properties}
\newtheorem{Sub}[Thm]{Sublemma}

\theoremstyle{definition}
\newtheorem{Def}[Thm]{Definition}
\newtheorem{Exm}[Thm]{Example}
\newtheorem{Exs}[Thm]{Examples}

\theoremstyle{remark}
\newtheorem{Rem}[Thm]{Remark}
\newtheorem{Rms}[Thm]{Remarks}
\newtheorem*{Com}{Commentary}

\newcommand{\myEmail}{piotr.niemiec@uj.edu.pl}
\newcommand{\myAddress}{\noindent{}Piotr Niemiec\\{}Jagiellonian University\\{}Institute of Mathematics\\{}
   ul. \L{}ojasiewicza 6\\{}30-348 Krak\'{o}w\\{}Poland}
\newcommand{\myData}{\author[P. Niemiec]{Piotr Niemiec}\address{\myAddress}\email{\myEmail}}

\newcommand{\NNN}{\mathbb{N}}
\newcommand{\RRR}{\mathbb{R}}

\newcommand{\CCc}{\CMcal{C}}
\newcommand{\FFf}{\CMcal{F}}

\newcommand{\WWw}{\CMcal{W}}

\newcommand{\CcC}{\EuScript{C}}

\newcommand{\KkK}{\EuScript{K}}

\newcommand{\WwW}{\EuScript{W}}

\newcommand{\Dd}{\mathfrak{D}}

\newcommand{\mM}{\mathfrak{m}}

\newcommand{\SECT}[1]{\section{#1}\renewcommand{\theequation}{\thesection-\arabic{equation}}\setcounter{equation}{0}}

\newcounter{help}
\newcommand{\ITE}[3]{\ifthenelse{#1}{#2}{#3}}\newcommand{\ITEE}[3]{\ITE{\equal{#1}{#2}}{#3}{}}

\newcommand{\Contr}{\operatorname{Contr}}
\newcommand{\Iso}{\operatorname{Iso}}\newcommand{\Prob}{\operatorname{Prob}}

\newcommand{\Fix}{\operatorname{Fix}}

\newcommand{\leqsl}{\leqslant}\newcommand{\geqsl}{\geqslant}

\newcommand{\varempty}{\varnothing}\newcommand{\dd}{\colon}
\newcommand{\dint}[1]{\,\textup{d} #1}
\newcommand{\quotient}[2]{\makebox[1ex][l]{\raisebox{1ex}{$#1$}}\makebox[1ex][l]{$\diagup$}\raisebox{-1ex}{$#2$}}

\newcommand{\nl}{${}$\newline}

\newcommand{\THM}[1]{Theorem~\textup{\ref{thm:#1}}}
\newcommand{\DEF}[1]{Definition~\textup{\ref{def:#1}}}\newcommand{\COR}[1]{Corollary~\textup{\ref{cor:#1}}}
\newcommand{\LEM}[1]{Lemma~\textup{\ref{lem:#1}}}\newcommand{\PRO}[1]{Proposition~\textup{\ref{pro:#1}}}

\newenvironment{thm}[1]{\begin{Thm}\label{thm:#1}}{\end{Thm}}\newenvironment{lem}[1]{\begin{Lem}\label{lem:#1}}{\end{Lem}}
\newenvironment{cor}[1]{\begin{Cor}\label{cor:#1}}{\end{Cor}}\newenvironment{pro}[1]{\begin{Pro}\label{pro:#1}}{\end{Pro}}
\newenvironment{dfn}[1]{\begin{Def}\label{def:#1}}{\end{Def}}
\newenvironment{exm}[1]{\begin{Exm}\label{exm:#1}}{\end{Exm}}
\newenvironment{rem}[1]{\begin{Rem}\label{rem:#1}}{\end{Rem}}

\newcommand{\bibITEM}[2]{\ITE{\equal{#2}{}}{\bibitem{#1} }{\bibitem[#2]{#1} }}
\newcommand{\BIB}[8]{
   \bibITEM{#1}{#8} #2, \textit{#3}, #4{} \textbf{#5} (#6), #7.}
\newcommand{\myBIB}[6][P. Niemiec]{#1, \textit{#2}, #3{}\ITE{\equal{#4}{}}{}{ \textbf{#4}} (#5), #6.}
\newcommand{\BIb}[6]{
   \bibITEM{#1}{#6} #2, \textit{#3}, #4, #5.}
\newcommand{\BiB}[9]{
   \bibITEM{#1}{#9} #2, \textit{#3}, #4{} \textit{#5}, #6, #7, #8.}

\newcommand{\myBAPP}[3][P. Niemiec]{
   #1, \textit{#2}, #3}

\newcommand{\jRN}[2][]{
   \ITEE{#2}{ActaM}{\ITE{\equal{#1}{+}}
      {Acta Mathematica}{Acta Math.}}
   \ITEE{#2}{AdvM}{\ITE{\equal{#1}{+}}
      {Advances in Mathematics}{Adv. in Math.}}
   \ITEE{#2}{ACS}{\ITE{\equal{#1}{+}}
      {Applied Categorical Structures}{Appl. Categor. Struct.}}
   \ITEE{#2}{ActaSM}{\ITE{\equal{#1}{+}}
      {Acta Scientiarum Mathematicarum}{Acta Sci. Math.}}
   \ITEE{#2}{AmJM}{\ITE{\equal{#1}{+}}
      {American Journal of Mathematics}{Amer. J. Math.}}
   \ITEE{#2}{AmMMon}{\ITE{\equal{#1}{+}}
      {American Mathematical Monthly}{Amer. Math. Mon.}}
   \ITEE{#2}{AnnSciEcNormSupT}{\ITE{\equal{#1}{+}}
      {Annales Scientifiques de l'\'{E}cole Normale Sup\'{e}rieure (3)}{Ann. Sci. \'{E}c. Norm. Sup\'{e}r. (3)}}
   \ITEE{#2}{AnnM}{\ITE{\equal{#1}{+}}
      {Annals of Mathematics}{Ann. Math.}}
   \ITEE{#2}{AnnProb}{\ITE{\equal{#1}{+}}
      {The Annals of Probability}{Ann. Probab.}}
   \ITEE{#2}{AnnPALog}{\ITE{\equal{#1}{+}}
      {Annals of Pure and Applied Logic}{Ann. Pure Appl. Logic}}
   \ITEE{#2}{ArchM}{\ITE{\equal{#1}{+}}
      {Archiv der Mathematik}{Arch. Math.}}
   \ITEE{#2}{AttiAccLincRendNat}{\ITE{\equal{#1}{+}}
      {Atti della Accademia Nazionale dei Lincei. Rendiconti. Classe di Scienze Fisiche, Matematiche e Naturali}
      {Atti Accad. Naz. Lincei Rend. Cl. Sci. Fis. Mat. Nat.}}
   \ITEE{#2}{BAMS}{\ITE{\equal{#1}{+}}
      {Bulletin of the American Mathematical Society}{Bull. Amer. Math. Soc.}}
   \ITEE{#2}{BLondMS}{\ITE{\equal{#1}{+}}
      {Bulletin of the London Mathematical Sociecy}{Bull. Lond. Math. Soc.}}
   \ITEE{#2}{BAPolSSSM}{\ITE{\equal{#1}{+}}
      {Bulletin de l'Acad\'{e}mie Polonaise des Sciences. S\'{e}rie des Sciences Math\'{e}matiques}
      {Bull. Acad. Pol. Sci. S\'{e}r. Sci. Math.}}
   \ITEE{#2}{BullSM}{\ITE{\equal{#1}{+}}
      {Bulletin des Sciences Math\'{e}matiques}{Bull. Sci. Math.}}
   \ITEE{#2}{BullPol}{\ITE{\equal{#1}{+}}
      {Bulletin of the Polish Academy of Sciences: Mathematics}{Bull. Polish Acad. Sci. Math.}}
   \ITEE{#2}{CanadJM}{\ITE{\equal{#1}{+}}
      {Canadian Journal Mathematics}{Canad. J. Math.}}
   \ITEE{#2}{CollectM}{\ITE{\equal{#1}{+}}
      {Collectanea Mathematica}{Collect. Math.}}
   \ITEE{#2}{CMUC}{\ITE{\equal{#1}{+}}
      {Commentationes Mathematicae Universitatis Carolinae}{Comment. Math. Univ. Carolin.}}
   \ITEE{#2}{CRParis}{\ITE{\equal{#1}{+}}
      {C. R. Paris}{C. R. Paris}}
   \ITEE{#2}{CRASParis}{\ITE{\equal{#1}{+}}
      {Comptes Rendus de l'Acad\'{e}mie des Sciences. Paris}{C. R. Acad. Sci. Paris}}
   \ITEE{#2}{CEurJM}{\ITE{\equal{#1}{+}}
      {Central European Journal of Mathematics}{Cent. Eur. J. Math.}}
   \ITEE{#2}{CMHelv}{\ITE{\equal{#1}{+}}
      {Commentarii Mathematici Helvetici}{Comment. Math. Helv.}}
   \ITEE{#2}{CollM}{\ITE{\equal{#1}{+}}
      {Colloquium Mathematicum}{Coll. Math.}}
   \ITEE{#2}{ComposM}{\ITE{\equal{#1}{+}}
      {Compositio Mathematica}{Compos. Math.}}
   \ITEE{#2}{CzMJ}{\ITE{\equal{#1}{+}}
      {Czechoslovak Mathematical Journal}{Czech. Math. J.}}
   \ITEE{#2}{DissM}{\ITE{\equal{#1}{+}}
      {Dissertationes Mathematicae}{Diss. Math.}}
   \ITEE{#2}{DANSSSR}{\ITE{\equal{#1}{+}}
      {Doklady Akademii Nauk SSSR}{Dokl. Akad. Nauk SSSR}}
   \ITEE{#2}{DukeMJ}{\ITE{\equal{#1}{+}}
      {Duke Mathematical Journal}{Duke Math. J.}}
   \ITEE{#2}{ELA}{\ITE{\equal{#1}{+}}
      {The Electronic Journal of Linear Algebra}{Electron. J. Linear Algebra}}
   \ITEE{#2}{ExtrM}{\ITE{\equal{#1}{+}}
      {Extracta Mathematicae}{Extracta Math.}}
   \ITEE{#2}{FM}{\ITE{\equal{#1}{+}}
      {Fundamenta Mathematicae}{Fund. Math.}}
   \ITEE{#2}{FAA}{\ITE{\equal{#1}{+}}
      {Functional Analysis and its Applications}{Funct. Anal. Appl.}}
   \ITEE{#2}{FunkAnalPril}{\ITE{\equal{#1}{+}}
      {Funktsional'ny\u{\i} Analiz i Ego Prilozheniya}{Funkts. Anal. Prilozh.}}
   \ITEE{#2}{GTopA}{\ITE{\equal{#1}{+}}
      {General Topology and its Applications}{General Topol. Appl.}}
   \ITEE{#2}{IllinoisJM}{\ITE{\equal{#1}{+}}
      {Illinois Journal of Mathematics}{Illinois J. Math.}}
   \ITEE{#2}{IndagMP}{\ITE{\equal{#1}{+}}
      {Indagationes Mathematicae (Proceedings)}{Indagationes Math. Proc.}}
   \ITEE{#2}{InHauEtSPM}{\ITE{\equal{#1}{+}}
      {Inst. Hautes \'{E}tudes Sci. Publ. Math.}{Inst. Hautes \'{E}tudes Sci. Publ. Math.}}
   \ITEE{#2}{IEOT}{\ITE{\equal{#1}{+}}
      {Integral Equations and Operator Theory}{Integr. Equ. Oper. Theory}}
   \ITEE{#2}{IsraelJM}{\ITE{\equal{#1}{+}}
      {Israel Journal of Mathematics}{Israel J. Math.}}
   \ITEE{#2}{JAusMSA}{\ITE{\equal{#1}{+}}
      {Journal of the Australian Mathematical Society. Series A}{J. Aust. Math. Soc. Ser. A}}
   \ITEE{#2}{JCA}{\ITE{\equal{#1}{+}}
      {Journal of Convex Analysis}{J. Convex Anal.}}
   \ITEE{#2}{JChinUST}{\ITE{\equal{#1}{+}}
      {J. China Univ. Sci. Tech.}{J. China Univ. Sci. Tech.}}
   \ITEE{#2}{JFA}{\ITE{\equal{#1}{+}}
      {Journal of Functional Analysis}{J. Funct. Anal.}}
   \ITEE{#2}{JMAnApp}{\ITE{\equal{#1}{+}}
      {J. Math. Anal. Appl.}{J. Math. Anal. Appl.}}
   \ITEE{#2}{JOT}{\ITE{\equal{#1}{+}}
      {Journal of Operator Theory}{J. Operator Theory}}
   \ITEE{#2}{KodaiMSemRep}{\ITE{\equal{#1}{+}}
      {Kodai Math. Sem. Rep.}{Kodai Math. Sem. Rep.}}
   \ITEE{#2}{LAA}{\ITE{\equal{#1}{+}}
      {Linear Algebra and its Applications}{Linear Algebra Appl.}}
   \ITEE{#2}{LMLA}{\ITE{\equal{#1}{+}}
      {Linear and Multilinear Algebra}{Linear Multilinear Algebra}}
   \ITEE{#2}{MLQ}{\ITE{\equal{#1}{+}}
      {Mathematical Logic Quarterly}{Math. Log. Q.}}
   \ITEE{#2}{MProcCambPhS}{\ITE{\equal{#1}{+}}
      {Mathematical Proceedings of the Cambridge Philosophical Society}{Math. Proc. Cambridge Phil. Soc.}}
   \ITEE{#2}{MMag}{\ITE{\equal{#1}{+}}
      {Mathematics Magazine}{Math. Mag.}}
   \ITEE{#2}{MSb}{\ITE{\equal{#1}{+}}
      {Matematicheski\u{\i} Sbornik}{Mat. Sb.}}
   \ITEE{#2}{MStud}{\ITE{\equal{#1}{+}}
      {Matematychni Studi\"{\i}}{Mat. Stud.}}
   \ITEE{#2}{MScand}{\ITE{\equal{#1}{+}}
      {Mathematica Scandinavica}{Math. Scand.}}
   \ITEE{#2}{MAnn}{\ITE{\equal{#1}{+}}
      {Mathematische Annalen}{Math. Ann.}}
   \ITEE{#2}{MAMS}{\ITE{\equal{#1}{+}}
      {Memoirs of the American Mathematical Society}{Mem. Amer. Math. Soc.}}
   \ITEE{#2}{MichMJ}{\ITE{\equal{#1}{+}}
      {Michigan Mathematical Journal}{Mich. Math. J.}}
   \ITEE{#2}{MonatM}{\ITE{\equal{#1}{+}}
      {Monatshefte f\"{u}r Mathematik}{Mh. Math.}}
   \ITEE{#2}{NonlinA}{\ITE{\equal{#1}{+}}
      {Nonlinear Analysis: Theory, Methods \& Applications}{Nonlinear Anal.}}
   \ITEE{#2}{OpusM}{\ITE{\equal{#1}{+}}
      {Opuscula Mathematica}{Opuscula Math.}}
   \ITEE{#2}{PacJM}{\ITE{\equal{#1}{+}}
      {Pacific Journal of Mathematics}{Pacific J. Math.}}
   \ITEE{#2}{PeriodMHung}{\ITE{\equal{#1}{+}}
      {Periodica Mathematica Hungarica}{Period. Math. Hungarica}}
   \ITEE{#2}{PAMS}{\ITE{\equal{#1}{+}}
      {Proceedings of the American Mathematical Society}{Proc. Amer. Math. Soc.}}
   \ITEE{#2}{ProcCambPhS}{\ITE{\equal{#1}{+}}
      {Proceedings of the Cambridge Philosophical Society}{Proc. Cambridge Phil. Soc.}}
   \ITEE{#2}{ProcImpAcadTokyo}{\ITE{\equal{#1}{+}}
      {Proc. Imp. Acad. Tokyo}{Proc. Imp. Acad. Tokyo}}
   \ITEE{#2}{ProcKonink}{\ITE{\equal{#1}{+}}
      {Proceedings of the Koninklijke Nederlandse Akademie van Wetenschappen}{Nederl. Akad. Wetensch. Proc. Ser. A}}
   \ITEE{#2}{PLondMS}{\ITE{\equal{#1}{+}}
      {Proceedings of the London Mathematical Society}{Proc. London Math. Soc.}}
   \ITEE{#2}{PNatlUSA}{\ITE{\equal{#1}{+}}
      {Proceedings of the National Academy of Sciences of the United States of America}{Proc. Natl. Acad. Sci. USA}}
   \ITEE{#2}{PublRIMSKyoto}{\ITE{\equal{#1}{+}}
      {Publ. Res. Inst. Math. Sci. Kyoto Univ.}{Publ. Res. Inst. Math. Sci.}}
   \ITEE{#2}{PWN}{\ITE{\equal{#1}{+}}
      {PWN -- Polish Scientific Publishers, Warszawa}{PWN -- Polish Scientific Publishers, Warszawa}}
   \ITEE{#2}{RCMP}{\ITE{\equal{#1}{+}}
      {Rendiconti del Circolo Matematico di Palermo}{Rend. Circ. Mat. Palermo}}
   \ITEE{#2}{RussMS}{\ITE{\equal{#1}{+}}
      {Russian Mathematical Surveys}{Russian Math. Surveys}}
   \ITEE{#2}{SbM}{\ITE{\equal{#1}{+}}
      {Sbornik: Mathematics}{Sb. Math.}}
   \ITEE{#2}{SciRepTokyoA}{\ITE{\equal{#1}{+}}
      {Science Reports of Tokyo Kyoiku Daigaku, Section A}{Sci. Rep. Tokyo Kyoiku Daigaku Sect. A}}
   \ITEE{#2}{SeminProbStras}{\ITE{\equal{#1}{+}}
      {S\'{e}minaire de probabilit\'{e}s de Strasbourg}{S\'{e}min. Prob. Strasbourg}}
   \ITEE{#2}{SIAMJMAA}{\ITE{\equal{#1}{+}}
      {SIAM Journal on Matrix Analysis and Applications}{SIAM J. Matrix Anal. Appl.}}
   \ITEE{#2}{SibirMZ}{\ITE{\equal{#1}{+}}
      {Sibirski\v{\i} Mat. \v{Z}hurnal}{Sibirsk. Mat. \v{Z}.}}
   \ITEE{#2}{SM}{\ITE{\equal{#1}{+}}
      {Studia Mathematica}{Studia Math.}}
   \ITEE{#2}{TAMS}{\ITE{\equal{#1}{+}}
      {Transactions of the American Mathematical Society}{Trans. Amer. Math. Soc.}}
   \ITEE{#2}{TohokuMJ}{\ITE{\equal{#1}{+}}
      {T\^{o}hoku Mathematical Journal}{T\^{o}hoku Math. J.}}
   \ITEE{#2}{TomskUnivRev}{\ITE{\equal{#1}{+}}
      {Tomsk Universitet Review}{Tomsk. Univ. Rev.}}
   \ITEE{#2}{TopA}{\ITE{\equal{#1}{+}}
      {Topology and its Applications}{Topology Appl.}}
   \ITEE{#2}{TsukubaJM}{\ITE{\equal{#1}{+}}
      {Tsukuba Journal of Mathematics}{Tsukuba J. Math.}}
   \ITEE{#2}{UspekhiMN}{\ITE{\equal{#1}{+}}
      {Uspekhi Matem. Nauk}{Uspekhi Mat. Nauk}}
   }

\newcommand{\paplist}[3][]{
   \ITEE{#3}{NIAkhiezer,IMGlazman1993}{
      \BIb{#2}{N.I. Akhiezer and I.M. Glazman}
         {Theory of Linear Operators in Hilbert Space}
         {Dover Publications, Inc., New York}{1993}{#1}}
   \ITEE{#3}{RDAnderson1967}{
      \BIB{#2}{R.D. Anderson}
         {On topological infinite deficiency}
         {\jRN{MichMJ}}{14}{1967}{365--383}{#1}}
   \ITEE{#3}{RDAnderson,JMcCharen1970}{
      \BIB{#2}{R.D. Anderson and J. McCharen}
         {On extending homeomorphisms to Fr\'{e}chet manifolds}
         {\jRN{PAMS}}{25}{1970}{283--289}{#1}}
   \ITEE{#3}{RDAnderson,DWCurtis,JVanMill1982}{
      \BIB{#2}{R.D. Anderson, D.W. Curtis, J. van Mill}
         {A fake topological Hilbert space}
         {\jRN{TAMS}}{272}{1982}{311--321}{#1}}
   \ITEE{#3}{RArens,JEells1956}{
      \BIB{#2}{R. Arens and J. Eells}
         {On embedding uniform and topological spaces}
         {\jRN{PacJM}}{6}{1956}{397--403}{#1}}
   \ITEE{#3}{NAronszajn,PPanitchpakdi1956}{
      \BIB{#2}{N. Aronszajn and P. Panitchpakdi}
         {Extension of uniformly continuous transformations and hyperconvex metric spaces}
         {\jRN{PacJM}}{6}{1956}{405--439}{#1}}
   \ITEE{#3}{KJBabenko1948}{
      \BIB{#2}{K.J. Babenko}
         {On conjugate functions}
         {\jRN{DANSSSR}}{62}{1948}{157--160}{#1}}
   \ITEE{#3}{TBanakh1995}{
      \BIB{#2}{T.O. Banakh}
         {Topology of spaces of probability measures, I}
         {\jRN{MStud}}{5}{1995}{65--87 (Russian)}{#1}}
   \ITEE{#3}{TBanakh1995a}{
      \BIB{#2}{T.O. Banakh}
         {Topology of spaces of probability measures, II}
         {\jRN{MStud}}{5}{1995}{88--106 (Russian)}{#1}}
   \ITEE{#3}{TBanakh1998}{
      \BIB{#2}{T. Banakh}
         {Characterization of spaces admitting a homotopy dense embedding into a Hilbert manifold}
         {\jRN{TopA}}{86}{1998}{123--131}{#1}}
   \ITEE{#3}{TBanakh,TNRadul1997}{
      \BIB{#2}{T.O. Banakh and T.N. Radul}
         {Topology of spaces of probability measures}
         {\jRN{SbM}}{188}{1997}{973--995}{#1}}
   \ITEE{#3}{TBanakh,TRadul,MZarichnyi1996}{
      \BIb{#2}{T. Banakh, T. Radul, M. Zarichnyi}
         {Absorbing sets in infinite-dimensional manifolds}
         {VNTL Publishers, Lviv}{1996}{#1}}
   \ITEE{#3}{TBanakh,IZarichnyy2008}{
      \BIB{#2}{T. Banakh and I. Zarichnyy}
         {Topological groups and convex sets homeomorphic to non-separable Hilbert spaces}
         {\jRN{CEurJM}}{6}{2008}{77--86}{#1}}
   \ITEE{#3}{HBecker,ASKechris1996}{
      \BIb{#2}{H. Becker and A.S. Kechris}
         {The Descriptive Set Theory of Polish Group Actions \textup{(London Math. Soc. Lecture Note Series, vol. 232)}}
         {University Press, Cambridge}{1996}{#1}}
   \ITEE{#3}{NEBenamara,NNikolski1999}{
      \BIB{#2}{N.E. Benamara and N. Nikolski}
         {Resolvent tests for similarity to a normal operator}
         {\jRN{PLondMS}}{78}{1999}{585--626}{#1}}
   \ITEE{#3}{YBenyamini,JLindenstrauss2000}{
      \BIb{#2}{Y. Benyamini and J. Lindenstrauss}
         {Geometric nonlinear functional analysis I}
         {AMS Colloquium Publications 48}{2000}{#1}}
   \ITEE{#3}{SKBerberian1974}{
      \BIb{#2}{S.K. Berberian}
         {Lectures in Functional Analysis and Operator Theory}
         {Graduate Texts in Mathematics 15, Springer-Verlag, New York}{1974}{#1}}
   \ITEE{#3}{SNBernstein1954}{
      \BIb{#2}{S.N. Bernstein}
         {Collected Works II}
         {Akad. Nauk SSSR, Moscow}{1954 (Russian)}{#1}}
   \ITEE{#3}{CzBessaga,APelczynski1972}{
      \BIB{#2}{Cz. Bessaga and A. Pe\l{}czy\'{n}ski}
         {On spaces of measurable functions}
         {\jRN{SM}}{44}{1972}{597--615}{#1}}
   \ITEE{#3}{CzBessaga,APelczynski1975}{
      \BIb{#2}{Cz. Bessaga and A. Pe\l{}czy\'{n}ski}
         {Selected topics in infinite-dimensional topology}
         {\jRN{PWN}}{1975}{#1}}
   \ITEE{#3}{MBestvina,JMogilski1986}{
      \BIB{#2}{M. Bestvina and J. Mogilski}
         {Characterizing certain incomplete infinite-dimensional absolute retracts}
         {\jRN{MichMJ}}{33}{1986}{291--313}{#1}}
   \ITEE{#3}{MBestvina,PBowers,JMogilsky,JWalsh1986}{
      \BIB{#2}{M. Bestvina, P. Bowers, J. Mogilsky, J. Walsh}
         {Characterization of Hilbert space manifolds revisited}
         {\jRN{TopA}}{24}{1986}{53--69}{#1}}
   \ITEE{#3}{RBhatia1997}{
      \BIb{#2}{R. Bhatia}
         {Matrix Analysis}
         {Springer, New York}{1997}{#1}}
   \ITEE{#3}{GBirkhoff1936}{
      \BIB{#2}{G. Birkhoff}
         {A note on topological groups}
         {\jRN{ComposM}}{3}{1936}{427--430}{#1}}
   \ITEE{#3}{MSBirman,MZSolomjak1987}{
      \BIb{#2}{M.S. Birman and M.Z. Solomjak}
         {Spectral Theory of Self-Adjoint Operators in Hilbert Space}
         {D. Reidel Publishing Co., Dordrecht}{1987}{#1}}
   \ITEE{#3}{EBishop1961}{
      \BIB{#2}{E. Bishop}
         {A generalization of the Stone-Weierstrass theorem}
         {\jRN{PacJM}}{11}{1961}{777--783}{#1}}
   \ITEE{#3}{BBlackadar2006}{
      \BIb{#2}{B. Blackadar}{Operator Algebras. Theory of $\CCc^*$-algebras and von Neumann algebras 
         \textup{(Encyclopaedia of Mathematical Sciences, vol. 122: Operator Algebras and Non-Commutative Geometry III)}}
         {Springer-Verlag, Berlin-Heidelberg}{2006}{#1}}
   \ITEE{#3}{JBlass,WHolsztynski1972}{
      \BIB{#2}{J. Blass and W. Holszty\'{n}ski}
         {Cubical polyhedra and homotopy III}
         {\jRN{AttiAccLincRendNat}}{53}{1972}{275--279}{#1}}
   \ITEE{#3}{FFBonsall,NJDuncan1973}{
      \BIb{#2}{F.F. Bonsall and N.J. Duncan}
         {Complete Normed Algebras}
         {Springer Verlag, Berlin}{1973}{#1}}
   \ITEE{#3}{NBourbaki2002}{
      \BIb{#2}{N. Bourbaki}
         {Lie Groups and Lie Algebras, Chapters 4--6}
         {Springer, New York}{2002}{#1}}
   \ITEE{#3}{PLBowers1989}{
      \BIB{#2}{P.L. Bowers}
         {Limitation topologies on function spaces}
         {\jRN{TAMS}}{314}{1989}{421--431}{#1}}
   \ITEE{#3}{AMBruckner,JBBruckner,BSThomson1997}{
      \BIb{#2}{A.M. Bruckner, J.B. Bruckner, B.S. Thomson}
         {Real Analysis}
         {Prentice-Hall, New Jersey}{1997}{#1}}
   \ITEE{#3}{PJCameron,AMVershik2006}{
      \BIB{#2}{P.J. Cameron and A.M. Vershik}
         {Some isometry groups of Urysohn space}
         {\jRN{AnnPALog}}{143}{2006}{70--78}{#1}}
   \ITEE{#3}{CCastaing1966}{
      \BIB{#2}{C. Castaing}
         {Quelques probl\`{e}mes de mesurabilit\'{e} li\'{e}es \`{a} la th\'{e}orie de la commande}
         {\jRN{CRParis}}{262}{1966}{409--411}{#1}}
   \ITEE{#3}{JAVanCasteren1980}{
      \BIB{#2}{J.A. van Casteren}
         {A problem of Sz.-Nagy}
         {\jRN{ActaSM}}{42}{1980}{189--194}{#1}}
   \ITEE{#3}{JAVanCasteren1983}{
      \BIB{#2}{J.A. van Casteren}
         {Operators similar to unitary or selfadjoint ones}
         {\jRN{PacJM}}{104}{1983}{241--255}{#1}}
   \ITEE{#3}{RCauty1994}{
      \BIB{#2}{R. Cauty}
         {Un espace m\'{e}trique lin\'{e}aire qui n'est pas un r\'{e}tracte absolu}
         {\jRN{FM}}{146}{1994}{85--99, (French)}{#1}}
   \ITEE{#3}{TAChapman1971}{
      \BIB{#2}{T.A. Chapman}
         {Deficiency in infinite-dimensional manifolds}
         {\jRN{GTopA}}{1}{1971}{263--272}{#1}}
   \ITEE{#3}{TAChapman1976}{
      \BIb{#2}{T.A. Chapman}
         {Lectures on Hilbert cube manifolds}
         {C.B.M.S. Regional Conference Series in Math. No 28, Amer. Math. Soc.}{1976}{#1}}
   \ITEE{#3}{RBChuaqui1977}{
      \BIB{#2}{R.B. Chuaqui}
         {Measures invariant under a group of transformations}
         {\jRN{PacJM}}{68}{1977}{313--329}{#1}}
   \ITEE{#3}{JBConway1985}{
      \BIb{#2}{J.B. Conway}
         {A Course in Functional Analysis}
         {Springer-Verlag, New York}{1985}{#1}}
   \ITEE{#3}{JBConway2000}{
      \BIb{#2}{J.B. Conway}
         {A Course in Operator Theory}
         {(Graduate Studies in Mathematics, vol. 21) Amer. Math. Soc., Providence}{2000}{#1}}
   \ITEE{#3}{GCorach,AMaestripieri,MMbekhta2009}{
      \BIB{#2}{G. Corach, A. Maestripieri, M. Mbekhta}
         {Metric and homogeneous structure of closed range operators}
         {\jRN{JOT}}{61}{2009}{171--190}{#1}}
   \ITEE{#3}{MJCowen,RGDouglas1978}{
      \BIB{#2}{M.J. Cowen and R.G. Douglas}
         {Complex geometry and operator theory}
         {\jRN{ActaM}}{141}{1978}{187--261}{#1}}
   \ITEE{#3}{DWCurtis1985}{
      \BIB{#2}{D.W. Curtis}
         {Boundary sets in the Hilbert cube}
         {\jRN{TopA}}{20}{1985}{201--221}{#1}}
   \ITEE{#3}{MMDay1958}{
      \BIb{#2}{M.M. Day}
         {Normed Linear Spaces}
         {Springer Verlag, Berlin}{1958}{#1}}
   \ITEE{#3}{CDellacherie1967}{
      \BIB{#2}{C. Dellacherie}
         {Un compl\'{e}ment au th\'{e}or\`{e}me de Weierstrass-Stone}
         {\jRN{SeminProbStras}}{1}{1967}{52--53}{#1}}
   \ITEE{#3}{JJDijkstra1987}{
      \BIB{#2}{J.J. Dijkstra}
         {Strong negligibility of $\sigma$-compacta does not characterize Hilbert space}
         {\jRN{PacJM}}{127}{1987}{19--30}{#1}}
   \ITEE{#3}{JJDijkstra1990}{
      \BIB{#2}{J.J. Dijkstra}
         {Characterizing Hilbert space topology in terms of strong negligibility}
         {\jRN{ComposM}}{75}{1990}{299--306}{#1}}
   \ITEE{#3}{TDobrowolski,WMarciszewski2002}{
      \BIB{#2}{T. Dobrowolski and W. Marciszewski}
         {Failure of the Factor Theorem for Borel pre-Hilbert spaces}
         {\jRN{FM}}{175}{2002}{53--68}{#1}}
   \ITEE{#3}{TDobrowolski,JMogilski1990}{
      \BiB{#2}{T. Dobrowolski and J. Mogilski}
         {Problems on Topological Classification of Incomplete Metric Spaces}{Chapter 25 in:}
         {Open Problems in Topology}{J. van Mill and G.M. Reed (eds.), North-Holland Amsterdam}{1990}{411--429}{#1}}
   \ITEE{#3}{TDobrowolski,HTorunczyk1981}{
      \BIB{#2}{T. Dobrowolski and H. Toru\'{n}czyk}
         {Separable complete ANR's admitting a group structure are Hilbert manifolds}
         {\jRN{TopA}}{12}{1981}{229--235}{#1}}
   \ITEE{#3}{RGDouglas1966}{
      \BIB{#2}{R.G. Douglas}
         {On majorization, factorization and range inclusion of operators in Hilbert space}
         {\jRN{PAMS}}{17}{1966}{413--416}{#1}}
   \ITEE{#3}{CHDowker1947}{
      \BIB{#2}{C.H. Dowker}
         {Mapping theorems for non-compact spaces}
         {\jRN{AmJM}}{69}{1947}{200--242}{#1}}
   \ITEE{#3}{CHDowker1952}{
      \BIB{#2}{C.H. Dowker}
         {Topology of metric complexes}
         {\jRN{AmJM}}{74}{1952}{555--577}{#1}}
   \ITEE{#3}{JDugundji1951}{
      \BIB{#2}{J. Dugundji}
         {An extension of Tietze's theorem}
         {\jRN{PacJM}}{1}{1951}{353--367}{#1}}
   \ITEE{#3}{JDugundji1958}{
      \BIB{#2}{J. Dugundji}
         {Absolute neighborhood retracts and local connectedness for arbitrary metric spaces}
         {\jRN{ComposM}}{13}{1958}{229--246}{#1}}
   \ITEE{#3}{JDugundji1965}{
      \BIB{#2}{J. Dugundji}
         {Locally equiconnected spaces and absolute neighborhood retracts}
         {\jRN{FM}}{57}{1965}{187--193}{#1}}
   \ITEE{#3}{NDunford,JTSchwartz1958}{
      \BIb{#2}{N. Dunford and J.T. Schwartz}
         {Linear Operators, part I}
         {Interscience Publishers, New York}{1958}{#1}}
   \ITEE{#3}{NDunford,JTSchwartz1963}{
      \BIb{#2}{N. Dunford and J.T. Schwartz}
         {Linear Operators, part II}
         {Interscience Publishers, New York}{1963}{#1}}
   \ITEE{#3}{NDunford,JTSchwartz1971}{
      \BIb{#2}{N. Dunford and J.T. Schwartz}
         {Linear Operators, part III}
         {Wiley-Interscience, New York}{1971}{#1}}
   \ITEE{#3}{MLEaton,MDPerlman1977}{
      \BIB{#2}{M.L. Eaton and M.D. Perlman}
         {Reflection groups, generalized Schur functions and the geometry of majorization}
         {\jRN{AnnProb}}{5}{1977}{829--860}{#1}}
   \ITEE{#3}{EGEffros1965}{
      \BIB{#2}{E.G. Effros}
         {The Borel space of von Neumann algebras on a separable Hilbert space}
         {\jRN{PacJM}}{15}{1965}{1153--1164}{#1}}
   \ITEE{#3}{EGEffros1966}{
      \BIB{#2}{E.G. Effros}
         {Global structure in von Neumann algebras}
         {\jRN{TAMS}}{121}{1966}{434--454}{#1}}
   \ITEE{#3}{REspinola,MAKhamsi2001}{
      \BiB{#2}{R. Espinola and M.A. Khamsi}
         {Introduction to hyperconvex spaces}{Chapter XIII in:}{Handbook of Metric Fixed Point Theory}
         {W.A. Kirk and B. Sims (editors), Kluwer Academic Publishers}{2001}{391--435}{#1}}
   \ITEE{#3}{PAFillmore,JPWilliams1971}{
      \BIB{#2}{P.A. Fillmore and J.P. Williams}
         {On operator ranges}
         {\jRN{AdvM}}{7}{1971}{254--281}{#1}}
   \ITEE{#3}{JEells,NHKuiper1969}{
      \BIB{#2}{J. Eells and N.H. Kuiper}
         {Homotopy negligible subsets in infinite-dimensional manifolds}
         {\jRN{ComposM}}{21}{1969}{151--161}{#1}}
   \ITEE{#3}{REngelking1977}{
      \BIb{#2}{R. Engelking}
         {General Topology}
         {\jRN{PWN}}{1977}{#1}}
   \ITEE{#3}{REngelking1978}{
      \BIb{#2}{R. Engelking}
         {Dimension Theory}
         {\jRN{PWN}}{1978}{#1}}
   \ITEE{#3}{REngelking1989}{
      \BIb{#2}{R. Engelking}
         {General Topology. Revised and completed edition \textup{(Sigma series in pure mathematics, vol. 6)}}
         {Heldermann Verlag, Berlin}{1989}{#1}}
   \ITEE{#3}{PErdos,RDMauldin1976}{
      \BIB{#2}{P. Erd\"{o}s and R.D. Mauldin}
         {The nonexistence of certain invariant measures}
         {\jRN{PAMS}}{59}{1976}{321--322}{#1}}
   \ITEE{#3}{JErnest1976}{
      \BIB{#2}{J. Ernest}
         {Charting the operator terrain}
         {\jRN{MAMS}}{171}{1976}{207 pp}{#1}}
   \ITEE{#3}{RHFox1943}{
      \BIB{#2}{R.H. Fox}
         {On fiber spaces, II}
         {\jRN{BAMS}}{49}{1943}{733--735}{#1}}
   \ITEE{#3}{NAFriedman1970}{
      \BIb{#2}{N.A. Friedman}
         {Introduction to ergodic theory}
         {Van Nostrand Reinhold Company}{1970}{#1}}
   \ITEE{#3}{SGao,ASKechris2003}{
      \BIB{#2}{S. Gao and A.S. Kechris}
         {On the classification of Polish metric spaces up to isometry}
         {\jRN{MAMS}}{161}{2003}{viii+78}{#1}}
   \ITEE{#3}{MIGarrido,FMontalvo1991}{
      \BIB{#2}{M.I. Garrido and F. Montalvo}
         {On some generalizations of the Kakutani-Stone and Stone-Weierstrass theorems}
         {\jRN{ExtrM}}{6}{1991}{156--159}{#1}}
   \ITEE{#3}{LGe,JShen2002}{
      \BIB{#2}{L. Ge and J. Shen}
         {Generator problem for certain property T factors}
         {\jRN{PNAS}}{99}{2002}{565--567}{#1}}
   \ITEE{#3}{IMGelfand,MANaimark1943}{
      \BIB{#2}{I.M. Gelfand and M.A. Naimark}
         {On the embedding of normed rings into the ring of operators in Hilbert space}
         {\jRN{MSb}}{12}{1943}{197--213}{#1}}
   \ITEE{#3}{FGesztesy,MMalamud,MMitrea,SNaboko2009}{
      \BIB{#2}{F. Gesztesy, M. Malamud, M. Mitrea, S. Naboko}
         {Generalized polar decompositions for closed operators in Hilbert spaces and some applications}
         {\jRN{IEOT}}{64}{2009}{83--113}{#1}}
   \ITEE{#3}{LGillman,MJerison1960}{
      \BIb{#2}{L. Gillman and M. Jerison}
         {Rings of continuous functions}
         {New York}{1960}{#1}}
   \ITEE{#3}{JGlimm1960}{
      \BIB{#2}{J. Glimm}
         {A Stone-Weierstrass theorem for $\CCc^*$-algebras}
         {\jRN{AnnM}}{72}{1960}{216--244}{#1}}
   \ITEE{#3}{GGodefroy,NJKalton2003}{
      \BIB{#2}{G. Godefroy and N.J. Kalton}
         {Lipschitz-free Banach spaces}
         {\jRN{SM}}{159}{2003}{121--141}{#1}}
   \ITEE{#3}{ICGohberg,MGKrein1967}{
      \BIB{#2}{I.C. Gohberg and M.G. Krein}
         {On a description of contraction operators similar to unitary ones}
         {\jRN{FunkAnalPril}}{1}{1967}{38--60}{#1}}
   \ITEE{#3}{ELGriffinJr1953}{
      \BIB{#2}{E.L. Griffin Jr.}
         {Some contributions to the theory of rings of operators}
         {\jRN{TAMS}}{75}{1953}{471--504}{#1}}
   \ITEE{#3}{ELGriffinJr1955}{
      \BIB{#2}{E.L. Griffin Jr.}
         {Some contributions to the theory of rings of operators II}
         {\jRN{TAMS}}{79}{1955}{389--400}{#1}}
   \ITEE{#3}{MGromov1981}{
      \BIB{#2}{M. Gromov}
         {Groups of polynomial growth and expanding maps}
         {\jRN{InHauEtSPM}}{53}{1981}{53--73}{#1}}
   \ITEE{#3}{MGromov1999}{
      \BIb{#2}{M. Gromov}
         {Metric Structures for Riemannian and Non-Riemannian Spaces}
         {Progress in Math. \textbf{152}, Birkh\"{a}user}{1999}{#1}}
   \ITEE{#3}{JDeGroot1956}{
      \BIB{#2}{J. de Groot}
         {Non-archimedean metrics in topology}
         {\jRN{PAMS}}{7}{1956}{948--953}{#1}}
   \ITEE{#3}{LCGrove,CTBenson1985}{
      \BIb{#2}{L.C. Grove and C.T. Benson}
         {Finite Reflection Group}
         {2nd ed., Springer-Verlag}{1985}{#1}}
   \ITEE{#3}{VIGurarii1966}{
      \BIB{#2}{V.I. Gurari\v{\i}}{Spaces of universal placement, isotropic spaces and a problem of Mazur 
         on rotations of Banach spaces \textup{(Russian)}}
         {\jRN{SibirMZ}}{7}{1966}{1002--1013}{#1}}
   \ITEE{#3}{HHahn1932}{
      \BIb{#2}{H. Hahn}
         {Reelle Funktionen I}
         {Leipzig}{1932}{#1}}
   \ITEE{#3}{PRHalmos1950}{
      \BIb{#2}{P.R. Halmos}
         {Measure theory}
         {Van Nostrand, New York}{1950}{#1}}
   \ITEE{#3}{PRHalmos1951}{
      \BIb{#2}{P.R. Halmos}
         {Introduction to Hilbert Space and the Theory of Spectral Multiplicity}
         {Chelsea Publishing Company, New York}{1951}{#1}}
   \ITEE{#3}{PRHalmos1956}{
      \BIb{#2}{P.R. Halmos}
         {Lectures on Ergodic Theory}
         {Publ. Math. Soc. Japan, Tokyo}{1956}{#1}}
   \ITEE{#3}{PRHalmos1982}{
      \BIb{#2}{P.R. Halmos}
         {A Hilbert Space Problem Book}
         {Springer-Verlag New York Inc.}{1982}{#1}}
   \ITEE{#3}{RWHansell1972}{
      \BIB{#2}{R.W. Hansell}
         {On the nonseparable theory of Borel and Souslin sets}
         {\jRN{BAMS}}{78}{1972}{236--241}{#1}}
   \ITEE{#3}{FHausdorff1930}{
      \BIB{#2}{F. Hausdorff}
         {Erweiterung einer Hom\"{o}omorphie}
         {\jRN{FM}}{16}{1930}{353--360}{#1}}
   \ITEE{#3}{FHausdorff1934}{
      \BIB{#2}{F. Hausdorff}
         {\"{U}ber innere Abbildungen}
         {\jRN{FM}}{23}{1934}{279--291}{#1}}
   \ITEE{#3}{FHausdorff1938}{
      \BIB{#2}{F. Hausdorff}
         {Erweiterung einer stetigen Abbildung}
         {\jRN{FM}}{30}{1938}{40--47}{#1}}
   \ITEE{#3}{DWHenderson1971}{
      \BIB{#2}{D.W. Henderson}
         {Corrections and extensions of two papers about infinite-dimensional manifolds}
         {\jRN{GTopA}}{1}{1971}{321--327}{#1}}
   \ITEE{#3}{DWHenderson1975}{
      \BIB{#2}{D.W. Henderson}
         {$Z$-sets in ANR's}
         {\jRN{TAMS}}{213}{1975}{205--216}{#1}}
   \ITEE{#3}{DWHenderson,RMSchori1970}{
      \BIB{#2}{D.W. Henderson and R.M. Schori}
         {Topological classification of infinite-dimensional manifolds by homotopy type}
         {\jRN{BAMS}}{76}{1970}{121--124}{#1}}
   \ITEE{#3}{DWHenderson,JEWest1970}{
      \BIB{#2}{D.W. Henderson and J.E. West}
         {Triangulated infinite-dimensional manifolds}
         {\jRN{BAMS}}{76}{1970}{655--660}{#1}}
   \ITEE{#3}{BHoffmann1979}{
      \BIB{#2}{B. Hoffmann}
         {A compact contractible topological group is trivial}
         {\jRN{ArchM}}{32}{1979}{585--587}{#1}}
   \ITEE{#3}{DHofmann2002}{
      \BIB{#2}{D. Hofmann}
         {On a generalization of the Stone-Weierstrass theorem}
         {\jRN{ACS}}{10}{2002}{569--592}{#1}}
   \ITEE{#3}{GHognas,AMukherjea1995}{
      \BIb{#2}{G. H\"ogn\"as and A. Mukherjea}
         {Probability Measures on Semigroups. Convolution Products, Random Walks, and Random Matrices}
         {Plenum Press, New York}{1995}{#1}}
   \ITEE{#3}{MRHolmes1992}{
      \BIB{#2}{M.R. Holmes}
         {The universal separable metric space of Urysohn and isometric embeddings thereof in Banach spaces}
         {\jRN{FM}}{140}{1992}{199--223}{#1}}
   \ITEE{#3}{MRHolmes2008}{
      \BIB{#2}{M.R. Holmes}
         {The Urysohn space embeds in Banach spaces in just one way}
         {\jRN{TopA}}{155}{2008}{1479--1482}{#1}}
   \ITEE{#3}{RRHolmes,TYTam1999}{
      \BIB{#2}{R.R. Holmes and T.Y. Tam}
         {Distance to the convex hull of an orbit under the action of a compact group}
         {\jRN{JAusMSA}}{66}{1999}{331--357}{#1}}
   \ITEE{#3}{RHorn,RMathias1990}{
      \BIB{#2}{R. Horn and R. Mathias}
         {Cauchy-Schwartz inequalities associated with positive semidefinite matrices}
         {\jRN{LAA}}{142}{1990}{63--82}{#1}}
   \ITEE{#3}{GEHuhunaisvili1955}{
      \BIB{#2}{G.E. Huhunai\v{s}vili}
         {On a property of Urysohn's universal metric space}
         {\jRN{DANSSSR}}{101}{1955}{607--610 (Russian)}{#1}}
   \ITEE{#3}{JEHumphreys1990}{
      \BIb{#2}{J.E. Humphreys}
         {Reflection Groups and Coxeter Groups}
         {Cambridge University Press}{1990}{#1}}
   \ITEE{#3}{JRIsbell1964}{
      \BIB{#2}{J.R. Isbell}
         {Six theorems about injective metric spaces}
         {\jRN{CMHelv}}{39}{1964}{65--76}{#1}}
   \ITEE{#3}{SIzumino,YKato1985}{
      \BIB{#2}{S. Izumino and Y. Kato}
         {The closure of invertible operators on Hilbert space}
         {\jRN{ActaSM}}{49}{1985}{321--327}{#1}}
   \ITEE{#3}{CJiang2004}{
      \BIB{#2}{C. Jiang}
         {Similarity classification of Cowen-Douglas operators}
         {\jRN{CanadJM}}{56}{2004}{742--775}{#1}}
   \ITEE{#3}{WBJohnson,JLindenstrauss2001}{
      \BiB{#2}{W.B. Johnson and J. Lindenstrauss}{Basic Concepts in the Geometry of Banach Spaces}
         {Chapter 1 in:}{Handbook of the Geometry of Banach Spaces, Vol. 1}
         {W.B. Johnson and J. Lindenstrauss (editors), Elsevier Science B.V., Amsterdam}{2001}{1--84}{#1}}
   \ITEE{#3}{IBJung,JStochel2008}{
      \BIB{#2}{I.B. Jung and J. Stochel}
         {Subnormal operators whose adjoints have rich point spectrum}
         {\jRN{JFA}}{255}{2008}{1797--1816}{#1}}
   \ITEE{#3}{RVKadison,JRRingrose1983}{
      \BIb{#2}{R.V. Kadison and J.R. Ringrose}
         {Fundamentals of the Theory of Operator Algebras. Volume I: Elementary Theory}
         {Academic Press, Inc., New York-London}{1983}{#1}}
   \ITEE{#3}{RVKadison,JRRingrose1986}{
      \BIb{#2}{R.V. Kadison and J.R. Ringrose}
         {Fundamentals of the Theory of Operator Algebras. Volume II: Advanced Theory}
         {Academic Press, Inc., Orlando-London}{1986}{#1}}
   \ITEE{#3}{SKakutani1936}{
      \BIB{#2}{S. Kakutani}
         {\"{U}ber die Metrisation der topologischen Gruppen}
         {\jRN{ProcImpAcadTokyo}}{12}{1936}{82--84}{#1}}
   \ITEE{#3}{SKakutani1938}{
      \BIB{#2}{S. Kakutani}
         {Two fixed-point theorems concerning bicompact convex sets}
         {\jRN{ProcImpAcadTokyo}}{14}{1938}{242--245}{#1}}
   \ITEE{#3}{SKakutani1941}{
      \BIB{#2}{S. Kakutani}
         {Concrete representation of abstract L-spaces}
         {\jRN{AnnM}}{42}{1941}{523--537}{#1}}
   \ITEE{#3}{SKakutani1941a}{
      \BIB{#2}{S. Kakutani}
         {Concrete representation of abstract M-spaces}
         {\jRN{AnnM}}{42}{1941}{994--1024}{#1}}
   \ITEE{#3}{NKalton2007}{
      \BIB{#2}{N. Kalton}
         {Extending Lipschitz maps into $\CCc(K)$-spaces}
         {\jRN{IsraelJM}}{162}{2007}{275--315}{#1}}
   \ITEE{#3}{RKane2001}{
      \BIb{#2}{R. Kane}
         {Reflection Groups and Invariant Theory}
         {Canadian Mathematical Society, Springer}{2001}{#1}}
   \ITEE{#3}{VKannan,SRRaju1980}{
      \BIB{#2}{V. Kannan and S.R. Raju}
         {The nonexistence of invariant universal measures on semigroups}
         {\jRN{PAMS}}{78}{1980}{482--484}{#1}}
   \ITEE{#3}{IKaplansky1951}{
      \BIB{#2}{I. Kaplansky}
         {A theorem on rings of operators}
         {\jRN{PacJM}}{1}{1951}{227--232}{#1}}
   \ITEE{#3}{MKatetov1988}{
      \BiB{#2}{M. Kat\v{e}tov}{On universal metric spaces}{in: Frolik (ed.),}
         {General Topology and its Relations to Modern Analysis and Algebra VI. Proceedings of the Sixth Prague 
         Topological Symposium 1986}{Heldermann Verlag Berlin}{1988}{323--330}{#1}}
   \ITEE{#3}{YKatznelson1960}{
      \BIB{#2}{Y. Katznelson}
         {Sur les alg\'{e}bres dont les \'{e}l\'{e}ments non n\'{e}gatifs admettent des racines carr\'{e}es}
         {\jRN{AnnSciEcNormSupT}}{77}{1960}{167--174}{#1}}
   \ITEE{#3}{OHKeller1931}{
      \BIB{#2}{O.H. Keller}
         {Die Homoiomorphie der kompakten konvexen Mengen in Hilbertschen Raum}
         {\jRN{MAnn}}{105}{1931}{748--758}{#1}}
   \ITEE{#3}{MAKhamsi,WAKirk,CMartinez2000}{
      \BIB{#2}{M.A. Khamsi, W.A. Kirk, C. Martinez}
         {Fixed point and selection theorems in hyperconvex spaces}
         {\jRN{PAMS}}{128}{2000}{3275--3283}{#1}}
   \ITEE{#3}{ABKhararazishvili1998}{
      \BIb{#2}{A.B. Khararazishvili}
         {Transformation groups and invariant measures. Set-theoretic aspects}
         {World Scientific Publishing Co., Inc., River Edge, NJ}{1998}{#1}}
   \ITEE{#3}{YKijima1987}{
      \BIB{#2}{Y. Kijima}
         {Fixed points of nonexpansive self-maps of a compact metric space}
         {\jRN{JMAnApp}}{123}{1987}{114--116}{#1}}
   \ITEE{#3}{JKindler1995}{
      \BIB{#2}{J. Kindler}
         {Minimax theorems with applications to convex metric spaces}
         {\jRN{CollM}}{68}{1995}{179--186}{#1}}
   \ITEE{#3}{WAKirk1998}{
      \BIB{#2}{W.A. Kirk}
         {Hyperconvexity of $\RRR$-trees}
         {\jRN{FM}}{156}{1998}{67--72}{#1}}
   \ITEE{#3}{VLKleeJr1952}{
      \BIB{#2}{V.L. Klee Jr.}
         {Invariant metrics in groups (solution of a problem of Banach)}
         {\jRN{PAMS}}{3}{1952}{484--487}{#1}}
   \ITEE{#3}{HJKowalsky1957}{
      \BIB{#2}{H.J. Kowalsky}
         {Einbettung metrischer R\"{a}ume}
         {\jRN{ArchM}}{8}{1957}{336--339}{#1}}
   \ITEE{#3}{WKubis,MRubin2010}{
      \BIB{#2}{W. Kubi\'{s} and M. Rubin}
         {Extension and reconstruction theorems for the Urysohn universal metric space}
         {\jRN{CzMJ}}{60}{2010}{1--29}{#1}}
   \ITEE{#3}{KKuratowski1966}{
      \BIb{#2}{K. Kuratowski}
         {Topology. \textup{Vol. I}}
         {\jRN{PWN}}{1966}{#1}}
   \ITEE{#3}{KKuratowski,BKnaster1927}{
      \BIB{#2}{K. Kuratowski and B. Knaster}
         {A connected and connected im kleinen point set which contains no perfect subset}
         {\jRN{BAMS}}{33}{1927}{106--109}{#1}}
   \ITEE{#3}{KKuratowski,AMostowski1976}{
      \BIb{#2}{K. Kuratowski and A. Mostowski}
         {Set Theory with an Introduction to Descriptive Set Theory}
         {\jRN{PWN}}{1976}{#1}}
   \ITEE{#3}{GLewicki1992}{
      \BIB{#2}{G. Lewicki}
         {Bernstein's ``lethargy'' theorem in metrizable topological linear spaces}
         {\jRN{MonatM}}{113}{1992}{213--226}{#1}}
   \ITEE{#3}{ASLewis1996}{
      \BIB{#2}{A.S. Lewis}
         {Group invariance and convex matrix analysis}
         {\jRN{SIAMJMAA}}{17}{1996}{927--949}{#1}}
   \ITEE{#3}{C-KLi,N-KTsing1991}{
      \BIB{#2}{C.-K. Li and N.-K. Tsing}
         {$G$-invariant norms and $G(c)$-radii}
         {\jRN{LAA}}{150}{1991}{179--194}{#1}}
   \ITEE{#3}{AJLazar,JLindenstrauss1971}{
      \BIB{#2}{A.J. Lazar and J. Lindenstrauss}
         {Banach spaces whose duals are $L_1$ spaces and their representing matrices}
         {\jRN{ActaM}}{126}{1971}{165--193}{#1}}
   \ITEE{#3}{EHLieb,MLoss1997}{
      \BIb{#2}{E.H. Lieb and M. Loss}
         {Analysis \textup{(Graduate Studies in Mathematics, vol. 14)}}
         {Amer. Math. Soc., Providence, RI}{1997}{#1}}
   \ITEE{#3}{ALindenbaum1926}{
      \BIB{#2}{A. Lindenbaum}
         {Contributions \`{a} l'\'{e}tude de l'espace m\'{e}trique I}
         {\jRN{FM}}{8}{1926}{209--222}{#1}}
   \ITEE{#3}{DLindenstrauss,LTzafriri1971}{
      \BIB{#2}{D. Lindenstrauss and L. Tzafriri}
         {On the complemented subspaces problem}
         {\jRN{IsraelJM}}{9}{1971}{263--269}{#1}}
   \ITEE{#3}{LHLoomis1945}{
      \BIB{#2}{L.H. Loomis}
         {Abstract congruence and the uniqueness of Haar measure}
         {\jRN{AnnM}}{46}{1945}{348--355}{#1}}
   \ITEE{#3}{LHLoomis1949}{
      \BIB{#2}{L.H. Loomis}
         {Haar measure in uniform structures}
         {\jRN{DukeMJ}}{16}{1949}{193--208}{#1}}
   \ITEE{#3}{ERLorch1939}{
      \BIB{#2}{E.R. Lorch}
         {Bicontinuous linear transformation in certain vector spaces}
         {\jRN{BAMS}}{45}{1939}{564--569}{#1}}
   \ITEE{#3}{ATLundell,SWeingram1969}{
      \BIb{#2}{A.T. Lundell and S. Weingram}
         {The topology of CW-complexes}
         {Litton Educ. Publ.}{1969}{#1}}
   \ITEE{#3}{WLusky1976}{
      \BIB{#2}{W. Lusky}
         {The Gurarij spaces are unique}
         {\jRN{ArchM}}{27}{1976}{627--635}{#1}}
   \ITEE{#3}{WLusky1977}{
      \BIB{#2}{W. Lusky}
         {On separable Lindenstrauss spaces}
         {\jRN{JFA}}{26}{1977}{103--120}{#1}}
   \ITEE{#3}{DMaharam1942}{
      \BIB{#2}{D. Maharam}
         {On homogeneous measure algebras}
         {\jRN{PNatlUSA}}{28}{1942}{108--111}{#1}}
   \ITEE{#3}{MMalicki,SSolecki2009}{
      \BIB{#2}{M. Malicki and S. Solecki}
         {Isometry groups of separable metric spaces}
         {\jRN{MProcCambPhS}}{146}{2009}{67--81}{#1}}
   \ITEE{#3}{PMankiewicz1972}{
      \BIB{#2}{P. Mankiewicz}
         {On extension of isometries in normed linear spaces}
         {\jRN{BAPolSSSM}}{20}{1972}{367--371}{#1}}
   \ITEE{#3}{JMartinezMaurica,MTPellon1987}{
      \BIB{#2}{J. Martinez-Maurica and M.T. Pell\'{o}n}
         {Non-archimedean Chebyshev centers}
         {\jRN{IndagMP}}{90}{1987}{417--421}{#1}}
   \ITEE{#3}{KMaurin1980}{
      \BIb{#2}{K. Maurin}
         {Analysis, Part II}
         {D. Reidel, Dordrecht-Boston-London}{1980}{#1}}
   \ITEE{#3}{SMazur,SUlam1932}{
      \BIB{#2}{S. Mazur and S. Ulam}
         {Sur les transformationes isom\'{e}triques d'espaces vectoriels norm\'{e}s}
         {\jRN{CRASParis}}{194}{1932}{946--948}{#1}}
   \ITEE{#3}{SMazurkiewicz1920}{
      \BIB{#2}{S. Mazurkiewicz}
         {Sur les lignes de Jordan}
         {\jRN{FM}}{1}{1920}{166--209}{#1}}
   \ITEE{#3}{SMazurkiewicz,WSierpinski1920}{
      \BIB{#2}{S. Mazurkiewicz and W. Sierpi\'{n}ski}
         {Contributions a la topologie des ensembles denombrables}
         {\jRN{FM}}{1}{1920}{17--27}{#1}}
   \ITEE{#3}{MMbekhta1992}{
      \BIB{#2}{M. Mbekhta}
         {Sur la structure des composantes connexes semi-Fredholm de $B(H)$}
         {\jRN{PAMS}}{116}{1992}{521--524}{#1}}
   \ITEE{#3}{JEMcCarthy1996}{
      \BIB{#2}{J.E. McCarthy}
         {Boundary values and Cowen-Douglas curvature}
         {\jRN{JFA}}{137}{1996}{1--18}{#1}}
   \ITEE{#3}{JMelleray2007}{
      \BIB{#2}{J. Melleray}
         {Computing the complexity of the relation of isometry between separable Banach spaces}
         {\jRN{MLQ}}{53}{2007}{128--131}{#1}}
   \ITEE{#3}{JMelleray2007a}{
      \BIB{#2}{J. Melleray}
         {On the geometry of Urysohn's universal metric space}
         {\jRN{TopA}}{154}{2007}{384--403}{#1}}
   \ITEE{#3}{JMelleray2008}{
      \BIB{#2}{J. Melleray}
         {Some geometric and dynamical properties of the Urysohn space}
         {\jRN{TopA}}{155}{2008}{1531--1560}{#1}}
   \ITEE{#3}{JMelleray,FVPetrov,AMVershik2008}{
      \BIB{#2}{J. Melleray, F.V. Petrov, A.M. Vershik}
         {Linearly rigid metric spaces and the embedding problem}
         {\jRN{FM}}{199}{2008}{177--194}{#1}}
   \ITEE{#3}{EMichael1953}{
      \BIB{#2}{E. Michael}
         {Some extension theorems for continuous functions}
         {\jRN{PacJM}}{3}{1953}{789--806}{#1}}
   \ITEE{#3}{EMichael1954}{
      \BIB{#2}{E. Michael}
         {Local properties of topological spaces}
         {\jRN{DukeMJ}}{21}{1954}{163--171}{#1}}
   \ITEE{#3}{EMichael1956}{
      \BIB{#2}{E. Michael}
         {Selected selection theorems}
         {\jRN{AmMMon}}{58}{1956}{233--238}{#1}}
   \ITEE{#3}{EMichael1956a}{
      \BIB{#2}{E. Michael}
         {Continuous selections. I}
         {\jRN{AnnM}}{63}{1956}{361--382}{#1}}
   \ITEE{#3}{EMichael1956b}{
      \BIB{#2}{E. Michael}
         {Continuous selections. II}
         {\jRN{AnnM}}{64}{1956}{562--580}{#1}}
   \ITEE{#3}{EMichael1959}{
      \BIB{#2}{E. Michael}
         {A theorem on semi-continuous set-valued functions}
         {\jRN{DukeMJ}}{26}{1959}{647--652}{#1}}
   \ITEE{#3}{JVanMill1986}{
      \BIB{#2}{J. van Mill}
         {Another counterexample in ANR theory}
         {\jRN{PAMS}}{97}{1986}{136--138}{#1}}
   \ITEE{#3}{JVanMill2001}{
      \BIb{#2}{J. van Mill}
         {The Infinite-Dimensional Topology of Function Spaces 
         \textup{(North-Holland Mathematical Library, vol. 64)}}
         {Elsevier, Amsterdam}{2001}{#1}}
   \ITEE{#3}{WMlak1991}{
      \BIb{#2}{W. Mlak}
         {Hilbert Spaces and Operator Theory}
         {PWN --- Polish Scientific Publishers and Kluwer Academic Publishers, Warszawa-Dordrecht}{1991}{#1}}
   \ITEE{#3}{JMogilski1979}{
      \BIB{#2}{J. Mogilski}
         {$CE$-decomposition of $l_2$-manifolds}
         {\jRN{BAPolSSSM}}{27}{1979}{309--314}{#1}}
   \ITEE{#3}{RLMoore1916}{
      \BIB{#2}{R.L. Moore}
         {On the foundations of plane analysis situs}
         {\jRN{TAMS}}{17}{1916}{131--164}{#1}}
   \ITEE{#3}{KMorita1955}{
      \BIB{#2}{K. Morita}
         {A condition for the metrizability of topological spaces and for $n$-dimensionality}
         {\jRN{SciRepTokyoA}}{5}{1955}{33--36}{#1}}
   \ITEE{#3}{AMukherjea,NATserpes1976}{
      \BIb{#2}{A. Mukherjea and N.A. Tserpes}
         {Measures on topological semigroups}
         {Springer Lecture Notes in Math. Vol. 547, Berlin}{1976}{#1}}
   \ITEE{#3}{JMycielski1974}{
      \BIB{#2}{J. Mycielski}
         {Remarks on invariant measures in metric spaces}
         {\jRN{CollM}}{32}{1974}{105--112}{#1}}
   \ITEE{#3}{SNNaboko1984}{
      \BIB{#2}{S.N. Naboko}
         {Conditions for similarity to unitary and selfadjoint operators}
         {\jRN{FunkAnalPril}}{18}{1984}{16--27}{#1}}
   \ITEE{#3}{LNachbin1965}{
      \BIb{#2}{L. Nachbin}
         {The Haar Integral}
         {D. Van Nostrand Company, Inc., Princeton-New Jersey-Toronto-New York-London}{1965}{#1}}
   \ITEE{#3}{TDNarang,SKGarg1991}{
      \BIB{#2}{T.D. Narang and S.K. Garg}
         {On the uniqueness of best approximation in non-archimedian spaces}
         {\jRN{PeriodMHung}}{22}{1991}{121--124}{#1}}
   \ITEE{#3}{JVonNeumann1930}{
      \BIB{#2}{J. von Neumann}
         {Zur Algebra der Funktionaloperationen und Theorie der normalen Operatoren}
         {\jRN{MAnn}}{102}{1930}{370--427}{#1}}
   \ITEE{#3}{JVonNeumann1934}{
      \BIB{#2}{J. von Neumann}
         {Zum Haarschen Mass in topologischen Gruppen}
         {\jRN{ComposM}}{1}{1934}{106--114}{#1}}
   \ITEE{#3}{JVonNeumann1937}{
      \BiB{#2}{J. von Neumann}
         {Some matrix-inequalities and metrization of matrix-space}{\jRN{TomskUnivRev}{} \textbf{1} (1937), 286--300; 
         in }{Collected Works}{Pergamon, New York}{1962}{Vol. 4, 205--219}{#1}}
   \ITEE{#3}{JVonNeumann1949}{
      \BIB{#2}{J. von Neumann}
         {On Rings of Operators. Reduction Theory}
         {\jRN{AnnM}}{50}{1949}{401--485}{#1}}
   \ITEE{#3}{pn2002}{\bibITEM{#2}{#1} \mypaplist{pn1}}
   \ITEE{#3}{pn2006a}{\bibITEM{#2}{#1} \mypaplist{pn2}}
   \ITEE{#3}{pn2006b}{\bibITEM{#2}{#1} \mypaplist{pn3}}
   \ITEE{#3}{pn2007}{\bibITEM{#2}{#1} \mypaplist{pn4}}
   \ITEE{#3}{pn2008a}{\bibITEM{#2}{#1} \mypaplist{pn5}}
   \ITEE{#3}{pn2008b}{\bibITEM{#2}{#1} \mypaplist{pn6}}
   \ITEE{#3}{pn2009a}{\bibITEM{#2}{#1} \mypaplist{pn7}}
   \ITEE{#3}{pn2009b}{\bibITEM{#2}{#1} \mypaplist{pn8}}
   \ITEE{#3}{pn2009c}{\bibITEM{#2}{#1} \mypaplist{pn9}}
   \ITEE{#3}{pn2010a}{\bibITEM{#2}{#1} \mypaplist{pn12}}
   \ITEE{#3}{pn2010b}{\bibITEM{#2}{#1} \mypaplist{pn13}}
   \ITEE{#3}{pn2011a}{\bibITEM{#2}{#1} \mypaplist{pn10}}
   \ITEE{#3}{pn2011b}{\bibITEM{#2}{#1} \mypaplist{pn15}}
   \ITEE{#3}{pn2011c}{\bibITEM{#2}{#1} \mypaplist{pn16}}
   \ITEE{#3}{pn2009x}{
      \bibITEM{#2}{#1} \mypaplist{pn11}}
   \ITEE{#3}{pn2010x}{
      \bibITEM{#2}{#1} \mypaplist{pn14}}
   \ITEE{#3}{pnXXXXa}{
      \bibITEM{#2}{#1} \mypaplist{pnX1}}
   \ITEE{#3}{pnXXXXb}{
      \bibITEM{#2}{#1} \mypaplist{pnX2}}
   \ITEE{#3}{pnXXXXc}{
      \bibITEM{#2}{#1} \mypaplist{pnX3}}
   \ITEE{#3}{pnXXXXd}{
      \bibITEM{#2}{#1} \mypaplist{pnX13}}
   \ITEE{#3}{MNiezgoda1998}{
      \BIB{#2}{M. Niezgoda}
         {Group majorization and Schur type inequalities}
         {\jRN{LAA}}{268}{1998}{9--30}{#1}}
   \ITEE{#3}{MNiezgoda1998a}{
      \BIB{#2}{M. Niezgoda}
         {An analytical characterization of effective and of irreducible groups inducing cone orderings}
         {\jRN{LAA}}{269}{1998}{105--114}{#1}}
   \ITEE{#3}{MNiezgoda,TYTam2001}{
      \BIB{#2}{M. Niezgoda and T.Y. Tam}
         {On norm property of $G(c)$-radii and Eaton triples}
         {\jRN{LAA}}{336}{2001}{119--130}{#1}}
   \ITEE{#3}{APazy1983}{
      \BIb{#2}{A. Pazy}{Semigroups of Linear Operators 
         and Applications to Partial Differential Equations \textup{(Applied Mathematical Sciences, vol. 44)}}
         {Springer-Verlag, New York}{1983}{#1}}
   \ITEE{#3}{APelc1982}{
      \BIB{#2}{A. Pelc}
         {Semiregular invariant measures on abelian groups}
         {\jRN{PAMS}}{86}{1982}{423--426}{#1}}
   \ITEE{#3}{RPenrose1955}{
      \BIB{#2}{R. Penrose}
         {A generalized inverse for matrices}
         {\jRN{ProcCambPhS}}{51}{1955}{406--413}{#1}}
   \ITEE{#3}{VPestov2006}{
      \BIb{#2}{V. Pestov}
         {Dynamics of infinite-dimensional groups. The Ramsey-Dvoretzky-Milman phenomenon}
         {University Lecture Series \textbf{40}, AMS, Providence, RI}{2006}{#1}}
   \ITEE{#3}{VPestov2007}{
      \BiB{#2}{V. Pestov}
         {Forty-plus annotated questions about large topological groups}
         {in:}{Open Problems in Topology II}{Elliot Pearl (editor), Elsevier B.V., Amsterdam}{2007}{439--450}{#1}}
   \ITEE{#3}{PVPetersen1993}{
      \BiB{#2}{P.V. Petersen}
         {Gromov-Hausdorff convergence of metric spaces}{in book:}{Differential Geometry: Riemannian Geometry 
         (Los Angeles, CA, 1990)}{Amer. Math. Soc., Providence, RI}{1993}{489--504}{#1}}
   \ITEE{#3}{DRamachandran,MMisiurewicz1982}{
      \BIB{#2}{D. Ramachandran and M. Misiurewicz}
         {Hopf's theorem on invariant measures for a group of transformations}
         {\jRN{SM}}{74}{1982}{183--189}{#1}}
   \ITEE{#3}{JMRosenblatt1974}{
      \BIB{#2}{J.M. Rosenblatt}
         {Equivalent invariant measures}
         {\jRN{IsraelJM}}{17}{1974}{261--270}{#1}}
   \ITEE{#3}{WRudin1962}{
      \BIb{#2}{W. Rudin}
         {Fourier Analysis on Groups \textup{(Interscience Tracts in Pure and Applied Mathematics, Number 12)}}
         {Interscience Publishers, New York}{1962}{#1}}
   \ITEE{#3}{WRudin1991}{
      \BIb{#2}{W. Rudin}
         {Functional Analysis}
         {McGraw-Hill Science}{1991}{#1}}
   \ITEE{#3}{TSaito1972}{
      \BiB{#2}{T. Sait\^{o}}{Generations of von Neumann algebras}
         {Lecture Notes in Math. vol. 247}{\textup{(}Lecture on Operator Algebras\textup{)}}
         {Springer, Berlin-Heidelberg-New York}{1972}{435--531}{#1}}
   \ITEE{#3}{KSakai,MYaguchi2003}{
      \BIB{#2}{K. Sakai and M. Yaguchi}
         {Characterizing manifolds modeled on certain dense subspaces of non-separable Hilbert spaces}
         {\jRN{TsukubaJM}}{27}{2003}{143--159}{#1}}
   \ITEE{#3}{SSakai1971}{
      \BIb{#2}{S. Sakai}
         {$\CCc^*$-Algebras and $\WWw^*$-Algebras}
         {Springer-Verlag, Berlin-Heidelberg-New York}{1971}{#1}}
   \ITEE{#3}{RSchori1971}{
      \BIB{#2}{R. Schori}
         {Topological stability for infinite-dimensional manifolds}
         {\jRN{ComposM}}{23}{1971}{87--100}{#1}}
   \ITEE{#3}{JTSchwartz1967}{
      \BIb{#2}{J.T. Schwartz}
         {$\WWw^*$-algebras}
         {Gordon and Breach, Science Publishers Inc., New York-London-Paris}{1967}{#1}}
   \ITEE{#3}{ZSemadeni1971}{
      \BIb{#2}{Z. Semadeni}
         {Banach Spaces of Continuous Functions (Vol. I)}
         {\jRN{PWN}}{1971}{#1}}
   \ITEE{#3}{JPSerre1951}{
      \BIB{#2}{J.-P. Serre}
         {Homologie singuli\`{e}re des espaces fibr\'{e}s}
         {\jRN{AnnM}}{54}{1951}{425--505}{#1}}
   \ITEE{#3}{DSherman2007}{
      \BIB{#2}{D. Sherman}
         {On the dimension theory of von Neumann algebras}
         {\jRN{MScand}}{101}{2007}{123--147}{#1}}
   \ITEE{#3}{WSierpinski1928}{
      \BIB{#2}{W. Sierpi\'{n}ski}
         {Sur les projections des ensembles compl\'{e}mentaires aux ensembles \textup{(A)}}
         {\jRN{FM}}{11}{1928}{117--122}{#1}}
   \ITEE{#3}{RCSteinlage1975}{
      \BIB{#2}{R.C. Steinlage}
         {On Haar measure in locally compact $T_2$ spaces}
         {\jRN{AmJM}}{97}{1975}{291--307}{#1}}
   \ITEE{#3}{JStochel,FHSzafraniec1989}{
      \BIB{#2}{J. Stochel and F.H. Szafraniec}
         {On normal extensions of unbounded operators. III. Spectral properties}
         {\jRN{PublRIMSKyoto}}{25}{1989}{105--139}{#1}}
   \ITEE{#3}{JStochel,FHSzafraniec1989a}{
      \BIB{#2}{J. Stochel and F.H. Szafraniec}
         {The normal part of an unbounded operator}
         {\jRN{ProcKonink}}{92}{1989}{495--503}{#1}}
   \ITEE{#3}{AHStone1962}{
      \BIB{#2}{A.H. Stone}
         {Absolute $\FFf_{\sigma}$-spaces}
         {\jRN{PAMS}}{13}{1962}{495--499}{#1}}
   \ITEE{#3}{AHStone1962a}{
      \BIB{#2}{A.H. Stone}
         {Non-separable Borel sets}
         {\jRN{DissM}}{28}{1962}{41 pages}{#1}}
   \ITEE{#3}{AHStone1972}{
      \BIB{#2}{A.H. Stone}
         {Non-separable Borel sets II}
         {\jRN{GTopA}}{2}{1972}{249--270}{#1}}
   \ITEE{#3}{MHStone1937}{
      \BIB{#2}{M.H. Stone}
         {Application of the theory of Boolean rings to general topology}
         {\jRN{TAMS}}{41}{1937}{375--481}{#1}}
   \ITEE{#3}{MHStone1948}{
      \BIB{#2}{M.H. Stone}
         {The generalized Weierstrass approximation theorem}
         {\jRN{MMag}}{21}{1948}{167--184}{#1}}
   \ITEE{#3}{BSz-Nagy1947}{
      \BIB{#2}{B. Sz.-Nagy}
         {On uniformly bounded linear transformations in Hilbert space}
         {\jRN{ActaSM}}{11}{1947}{152--157}{#1}}
   \ITEE{#3}{WTakahashi1970}{
      \BIB{#2}{W. Takahashi}
         {A convexity in metric space and nonexpansive mappings, I}
         {\jRN{KodaiMSemRep}}{22}{1970}{142--149}{#1}}
   \ITEE{#3}{MTakesaki2002}{
      \BIb{#2}{M. Takesaki}
         {Theory of Operator Algebras I \textup{(Encyclopaedia of Mathematical Sciences, Volume 124)}}
         {Springer-Verlag, Berlin-Heidelberg-New York}{2002}{#1}}
   \ITEE{#3}{MTakesaki2003}{
      \BIb{#2}{M. Takesaki}
         {Theory of Operator Algebras II \textup{(Encyclopaedia of Mathematical Sciences, Volume 125)}}
         {Springer-Verlag, Berlin-Heidelberg-New York}{2003}{#1}}
   \ITEE{#3}{MTakesaki2003a}{
      \BIb{#2}{M. Takesaki}
         {Theory of Operator Algebras III \textup{(Encyclopaedia of Mathematical Sciences, Volume 127)}}
         {Springer-Verlag, Berlin-Heidelberg-New York}{2003}{#1}}
   \ITEE{#3}{TYTam1999}{
      \BIB{#2}{T.Y. Tam}
         {An extension of a result of Lewis}
         {\jRN{ELA}}{5}{1999}{1--10}{#1}}
   \ITEE{#3}{TYTam2000}{
      \BIB{#2}{T.Y. Tam}
         {Group majorization, Eaton triples and numerical range}
         {\jRN{LMLA}}{47}{2000}{11--28}{#1}}
   \ITEE{#3}{TYTam2002}{
      \BIB{#2}{T.Y. Tam}
         {Generalized Schur-concave functions and Eaton triples}
         {\jRN{LMLA}}{50}{2002}{113--120}{#1}}
   \ITEE{#3}{TYTam,WCHill2001}{
      \BIB{#2}{T.Y. Tam and W.C. Hill}
         {On $G$-invariant norms}
         {\jRN{LAA}}{331}{2001}{101--112}{#1}}
   \ITEE{#3}{AFTiman,IAVestfrid1983}{
      \BIB{#2}{A.F. Timan and I.A. Vestfrid}
         {Any separable ultrametric space can be isometrically imbedded in $l_2$}
         {\jRN{FAA}}{17}{1983}{70--71}{#1}}
   \ITEE{#3}{JTomiyama1958}{
      \BIB{#2}{J. Tomiyama}
         {Generalized dimension function for $\WWw^*$-algebras of infinite type}
         {\jRN{TohokuMJ} (2)}{10}{1958}{121--129}{#1}}
   \ITEE{#3}{HTorunczyk1970}{
      \BIB{#2}{H. Toru\'{n}czyk}
         {Remarks on Anderson's paper ``On topological infinite deficiency''}
         {\jRN{FM}}{66}{1970}{393--401}{#1}}
   \ITEE{#3}{HTorunczyk1970a}{
      \BIb{#2}{H. Toru\'{n}czyk}
         {$G$-$K$-absorbing and skeletonized sets in metric spaces}
         {Ph.D. thesis, Inst. Math. Polish Acad. Sci., Warszawa}{1970}{#1}}
   \ITEE{#3}{HTorunczyk1972}{
      \BIB{#2}{H. Toru\'{n}czyk}
         {A short proof of Hausdorff's theorem on extending metrics}
         {\jRN{FM}}{77}{1972}{191--193}{#1}}
   \ITEE{#3}{HTorunczyk1974}{
      \BIB{#2}{H. Toru\'{n}czyk}
         {Absolute retracts as factors of normed linear spaces}
         {\jRN{FM}}{86}{1974}{53--67}{#1}}
   \ITEE{#3}{HTorunczyk1975}{
      \BIB{#2}{H. Toru\'{n}czyk}
         {On Cartesian factors and the topological classification of linear metric spaces}
         {\jRN{FM}}{88}{1975}{71--86}{#1}}
   \ITEE{#3}{HTorunczyk1978}{
      \BIB{#2}{H. Toru\'{n}czyk}
         {Concerning locally homotopy negligible sets and characterization of $l_2$-manifolds}
         {\jRN{FM}}{101}{1978}{93--110}{#1}}
   \ITEE{#3}{HTorunczyk1980}{
      \BiB{#2}{H. Toru\'{n}czyk}{Characterization of infinite-dimensional manifolds}{in:}
         {Proceedings of the International Conference on Geometric Topology (Warsaw, 1978)}
         {\jRN{PWN}}{1980}{431--437}{#1}}
   \ITEE{#3}{HTorunczyk1981}{
      \BIB{#2}{H. Toru\'{n}czyk}
         {Characterizing Hilbert space topology}
         {\jRN{FM}}{111}{1981}{247--262}{#1}}
   \ITEE{#3}{HTorunczyk1985}{
      \BIB{#2}{H. Toru\'{n}czyk}
         {A correction of two papers concerning Hilbert manifolds}
         {\jRN{FM}}{125}{1985}{89--93}{#1}}
   \ITEE{#3}{KTsuda1985}{
      \BIB{#2}{K. Tsuda}
         {A note on closed embeddings of finite dimensional metric spaces}
         {\jRN{BLondMS}}{17}{1985}{273--278}{#1}}
   \ITEE{#3}{PSUrysohn1925}{
      \BIB{#2}{P.S. Urysohn}
         {Sur un espace m\'{e}trique universel}
         {\jRN{CRASParis}}{180}{1925}{803--806}{#1}}
   \ITEE{#3}{PSUrysohn1927}{
      \BIB{#2}{P.S. Urysohn}
         {Sur un espace m\'{e}trique universel}
         {\jRN{BullSM}}{51}{1927}{43--64, 74--96}{#1}}
   \ITEE{#3}{VVUspenskij1986}{
      \BIB{#2}{V.V. Uspenskij}
         {A universal topological group with a countable basis}
         {\jRN{FAA}}{20}{1986}{86--87}{#1}}
   \ITEE{#3}{VVUspenskij1990}{
      \BIB{#2}{V.V. Uspenskij}
         {On the group of isometries of the Urysohn universal metric space}
         {\jRN{CMUC}}{31}{1990}{181--182}{#1}}
   \ITEE{#3}{VVUspenskij2004}{
      \BIB{#2}{V.V. Uspenskij}
         {The Urysohn universal metric space is homeomorphic to a Hilbert space}
         {\jRN{TopA}}{139}{2004}{145--149}{#1}}
   \ITEE{#3}{VVUspenskij2008}{
      \BIB{#2}{V.V. Uspenskij}
         {On subgroups of minimal topological groups}
         {\jRN{TopA}}{155}{2008}{1580--1606}{#1}}
   \ITEE{#3}{VSVaradarajan1963}{
      \BIB{#2}{V.S. Varadarajan}
         {Groups of automorphisms of Borel spaces}
         {\jRN{TAMS}}{109}{1963}{191--220}{#1}}
   \ITEE{#3}{AMVershik1998}{
      \BIB{#2}{A.M. Vershik}
         {The universal Urysohn space, Gromov's metric triples, and random metrics on the series of natural numbers}
         {\jRN{UspekhiMN}}{53}{1998}{57--64}{#1} English translation: \jRN{RussMS}{} \textbf{53} (1998), 921--928. 
         Correction: \jRN{UspekhiMN}{} \textbf{56} (2001), p. 207. English translation: \jRN{RussMS}{} \textbf{56} 
         (2001), p. 1015.}
   \ITEE{#3}{AMVershik2002}{
      \BIb{#2}{A.M. Vershik}
         {Random metric spaces and the universal Urysohn space}
         {Fundamental Mathematics Today. 10th anniversary of the Independent Moscow University. MCCME Publ.}{2002}{#1}}
   \ITEE{#3}{NWeaver1999}{
      \BIb{#2}{N. Weaver}
         {Lipschitz Algebras}
         {World Scientific}{1999}{#1}}
   \ITEE{#3}{JWeidmann1980}{
      \BIb{#2}{J. Weidmann}
         {Linear Operators in Hilbert Spaces}
         {(Graduate Texts in Mathematics, vol. 68) Springer-Verlag New York Inc.}{1980}{#1}}
   \ITEE{#3}{JEWest1969}{
      \BIB{#2}{J.E. West}
         {Approximating homotopies by isotopies in Fr\'{e}chet manifolds}
         {\jRN{BAMS}}{75}{1969}{1254--1257}{#1}}
   \ITEE{#3}{JEWest1969a}{
      \BIB{#2}{J.E. West}
         {Fixed-point sets of transformation groups on infinite-product spaces}
         {\jRN{PAMS}}{21}{1969}{575--582}{#1}}
   \ITEE{#3}{JEWest1970}{
      \BIB{#2}{J.E. West}
         {The ambient homeomorphy of infinite-dimensional Hilbert spaces}
         {\jRN{PacJM}}{34}{1970}{257--267}{#1}}
   \ITEE{#3}{JHCWhitehead1949}{
      \BIB{#2}{J.H.C. Whitehead}
         {Combinatorial homotopy I}
         {\jRN{BAMS}}{55}{1949}{213--245}{#1}}
   \ITEE{#3}{GTWhyburn1942}{
      \BIb{#2}{G. T. Whyburn}
         {Analytic Topology}
         {Amer. Math. Soc. Colloquium Publications (vol. XXVIII), New York}{1942}{#1}}
   \ITEE{#3}{WWogen1969}{
      \BIB{#2}{W. Wogen}
         {On generators for von Neumann algebras}
         {\jRN{BAMS}}{75}{1969}{95--99}{#1}}
   \ITEE{#3}{RYTWong1967}{
      \BIB{#2}{R.Y.T. Wong}
         {On homeomorphisms of certain infinite dimensional spaces}
         {\jRN{TAMS}}{128}{1967}{148--154}{#1}}
   \ITEE{#3}{LYang,JZhang1987}{
      \BIB{#2}{L. Yang and J. Zhang}
         {Average distance constants of some compact convex space}
         {\jRN{JChinUST}}{17}{1987}{17--23}{#1}}
   \ITEE{#3}{PZakrzewski1993}{
      \BIB{#2}{P. Zakrzewski}
         {The existence of invariant $\sigma$-finite measures for a group of transformations}
         {\jRN{IsraelJM}}{83}{1993}{275--287}{#1}}
   \ITEE{#3}{PZakrzewski2002}{
      \BIb{#2}{P. Zakrzewski}
         {Measures on Algebraic-Topological Structures, Handbook of Measure Thoery}
         {E. Pap, ed., Elsevier, Amsterdam}{2002, 1091--1130}{#1}}
   \ITEE{#3}{KZhu2000}{
      \BIB{#2}{K. Zhu}
         {Operators in Cowen-Douglas classes}
         {\jRN{IllinoisJM}}{44}{2000}{767--783}{#1}}
   }

\newcommand{\mypaplist}[2][]{
   \ITEE{#2}{pn1}{
      \myBIB{Separate and joint similarity to families of normal operators}
         {\jRN[#1]{SM}}{149}{2002}{39--62}}
   \ITEE{#2}{pn2}{
      \myBIB{Locally arcwise connected metrizable spaces with the fixed point property are complete-metrizable}
         {\jRN[#1]{TopA}}{153}{2006}{1639--1642}}
   \ITEE{#2}{pn3}{
      \myBIB{Invariant measures for equicontinuous semigroups of continuous transformations 
         of a compact Hausdorff space}{\jRN[#1]{TopA}}{153}{2006}{3373--3382}}
   \ITEE{#2}{pn4}{
      \myBIB{Approximation of the Hausdorff distance by the distance of continuous surjections}
         {\jRN[#1]{TopA}}{154}{2007}{655--664}}
   \ITEE{#2}{pn5}{
      \myBIB{Generalized Haar integral}
         {\jRN[#1]{TopA}}{155}{2008}{1323--1328}}
   \ITEE{#2}{pn6}{
      \myBIB{Integration and Lipschitz functions}
         {\jRN[#1]{RCMP}}{57}{2008}{391--399}}
   \ITEE{#2}{pn7}{
      \myBIB{Canonical Banach function spaces generated by Urysohn universal spaces. Measures as Lipschitz maps}
         {\jRN[#1]{SM}}{192}{2009}{97--110}}
   \ITEE{#2}{pn8}{
      \myBIB{Urysohn universal spaces as metric groups of exponent $2$}
         {\jRN[#1]{FM}}{204}{2009}{1--6}}
   \ITEE{#2}{pn9}{
      \myBIB{Central subsets of Urysohn universal spaces}
         {\jRN[#1]{CMUC}}{50}{2009}{445--461}}
   \ITEE{#2}{pn10}{
      \myBIB[P. Niemiec and T.Y. Tam]{A representation of $G$-in\-variant norms for Eaton triple}
         {\jRN[#1]{JCA}}{18}{2011}{59--65}}
   \ITEE{#2}{pn11}{
      \myBIB{Functor of extension of contractions on Urysohn universal spaces}
         {\jRN[#1]{ACS}}{}{2009}{\texttt{DOI: 10.1007/s10485-009-9218-z}}}
   \ITEE{#2}{pn12}{
      \myBIB{Ultra-$\mM$-separability}
         {\jRN[#1]{TopA}}{157}{2010}{669--673}}
   \ITEE{#2}{pn13}{
      \myBIB{Functor of extension of $\Lambda$-isometric maps between central subsets 
         of the unbounded Urysohn universal space}{\jRN[#1]{CMUC}}{51}{2010}{541--549}}
   \ITEE{#2}{pn14}{
      \myBIB{Normed topological pseudovector groups}{\jRN[#1]{ACS}}{}{2010}
         {\ITE{\equal{#1}{}}{\texttt{DOI: 10.1007/s10485\-010-9239-7}}{\texttt{DOI: 10.1007/s10485-010-9239-7}}}}
   \ITEE{#2}{pn15}{
      \myBIB{Topological structure of Urysohn universal spaces}
         {\jRN[#1]{TopA}}{158}{2011}{352--359}}
   \ITEE{#2}{pn16}{
      \myBIB{A note on invariant measures}
         {\jRN[#1]{OpusM}}{31}{2011}{425--431}}
   \ITEE{#2}{pnX1}{
      \myBAPP{Strengthened Stone-Weierstrass type theorem}
         {accepted for publication in \jRN[#1]{OpusM}}}
   \ITEE{#2}{pnX2}{
      \myBAPP{Functor of continuation in Hilbert cube and Hilbert space}
         {to appear in \jRN[#1]{FM}}}
   \ITEE{#2}{pnX3}{
      \myBAPP{Norm closures of orbits of bounded operators}
         {to appear.}}
   \ITEE{#2}{pnX6}{
      \myBAPP{Extending maps by injective $\sigma$-$Z$-maps in Hilbert manifolds}
         {to appear in \jRN[#1]{BullPol}}}
   \ITEE{#2}{pnX7}{
      \myBAPP{Spaces of measurable functions}
         {submitted to \jRN[#1]{CollectM}}}
   \ITEE{#2}{pnX8}{
      \myBAPP{Normal systems over ANR's, rigid embeddings and nonseparable absorbing sets}
         {submitted to \jRN[#1]{}}}
   \ITEE{#2}{pnX9}{
      \myBAPP{Borel structure of the spectrum of a closed operator}
         {submitted to \jRN[#1]{JFA}}}
   \ITEE{#2}{pnX10}{
      \myBAPP{Central points and measures and dense subsets of compact metric spaces}
         {submitted to \jRN[#1]{NonlinA}}}
   \ITEE{#2}{pnX11}{
      \myBAPP{Generalized absolute values and polar decompositions of a bounded operator}
         {submitted to \jRN[#1]{IEOT}.}}
   \ITEE{#2}{pnX12}{
      \myBAPP{Ultrametrics, extending of Lipschitz maps and nonexpansive selections}
         {submitted to \jRN[#1]{HJM}}}
   \ITEE{#2}{pnX13}{
      \myBAPP{A note on ANR's}
         {submitted to \jRN[#1]{TopA}}}
   \ITEE{#2}{pnX14}{
      \myBAPP{Problem with almost everywhere equality}
         {submitted to \jRN[#1]{ArchM}}}
   \ITEE{#2}{pnX15}{
      \myBAPP{Universal valued Abelian groups}
         {submitted to \jRN[#1]{MProcCambPhS}}}
   }

\begin{document}

\title[Central points and measures and dense subsets]
   {Central points and measures and dense subsets of compact metric spaces}
\myData
\begin{abstract}
For every nonempty compact convex subset $K$ of a normed linear space a (unique) point $c_K \in K$, called the generalized
Chebyshev center, is distinguished. It is shown that $c_K$ is a common fixed point for the isometry group of the metric
space $K$. With use of the generalized Chebyshev centers, the central measure $\mu_X$ of an arbitrary compact metric space
$X$ is defined. For a large class of compact metric spaces, including the interval $[0,1]$ and all compact metric groups,
another `central' measure is distinguished, which turns out to coincide with the Lebesgue measure and the Haar one
for the interval and a compact metric group, respectively. An idea of distinguishing infinitely many points forming
a dense subset of an arbitrary compact metric space is also presented.
\end{abstract}
\maketitle


\SECT{Introduction}

Distinguishing points, subsets or other `ingredients' related to spaces is important in many parts of mathematics,
including algebraic topology (homotopy groups), theory of Lipschitz functions (the base point), theory of locally compact
groups (the Haar measure, unique up to a constant factor). In most of algebraic structures the neutral element is
a naturally distinguished point. In other areas of mathematics distinguishing appears as a useful tool. For example,
the well-known \textit{Chebyshev center} of a nonempty compact convex subset of a \textbf{strictly} convex normed linear
space (i.e. such a space in which the unit sphere contains no segments \cite[page~30]{j-l}) finds an application in fixed
point theory, being a common fixed point for the isometry group of the convex set. The characteristic and important feature
of some of the above examples is the uniqueness, in a categorical or weaker sense, of the distinguished ingredients.
In such cases this distinguished ingredient may be seen as an integral part of the space (e.g. the Haar measure
of a locally compact group or the neutral element of an algebraic structure), while in the others it plays an additional
role (e.g. in homotopy groups, spaces of Lipschitz functions). In the latter cases the distinguishing is just a necessity
and it hardly ever finds applications. The foregoing examples show that the situation changes when the distinguished
ingredient turns out to be uniquely determined by some natural conditions. (The precise meaning of this in category
of metric spaces shall be explained in the next section.)\par
The aim of the recent paper is to present a few results dealing with constructive `applied' distinguishings. In particular,
we shall show that every nonempty compact metric space $X$ isometric to a convex subset of a normed linear space (even
of a more general class, containing all metric $\RRR$-trees) contains a unique point $c_X$ (called the \textit{generalized
Chebyshev center}) which is in a sense its center. As an application of this, we shall prove that the isometry group
of each such space has a common fixed point. This gives a constructive proof of Kakutani's fixed point theorem in a special
case. Details are included in Section~3.\par
In Section 4 we shall apply the results of the previous part to an arbitrary (nonempty) compact metric space $X$ in order
to define the \textit{central} (probability Borel) measure $\mu_X$ of $X$ by means of the so-called Kantorovich
(or Kantorovicz-Rubenstein, cf. \cite{weaver}) metric induced by the metric of $X$. In case of a compact metric group $G$,
$\mu_G$ turns out to be the Haar measure of $G$ and thus we shall obtain an alternative proof of the Haar measure theorem
for compact metrizable groups. However, the problem of whether $\mu_{[0,1]}$ is the one-dimensional Lebesgue measure
we leave as open. Section~5 deals with the so-called \textit{quasi-nilpotent} compact metric spaces for which we shall
prove another result on distinguishing measures. As a special case we shall obtain the characterizations of the Lebesgue
measure on $[0,1]$ and (again) the Haar measure of a compact metric group. The last, sixth, part is devoted
to distinguishing countable dense subsets in arbitrary compact metric spaces, which is related to theory of random metric
spaces (see e.g. \cite{vershik1,vershik2}).

\SECT{Preliminaries}

In this paper we deal with categories of metric spaces with additional structures in which every \textit{isomorphism}
between spaces is an isometric function between them. For simplicity, let us call each such a category
an \textit{iso-category}. We shall write $K \in \KkK$ to express that $K$ is a metric space with an additional structure
which belongs to an iso-category $\KkK$.\par
Let $\KkK$ be an iso-category. For any two members $X$ and $Y$ of $\KkK$ let $\Iso_{\KkK}(X,Y)$ stand for the set of all
isomorphisms of $X$ onto $Y$. We write $\Iso_{\KkK}(X)$ for $\Iso_{\KkK}(X,X)$. If no additional structures on metric
spaces are needed to describe the category $\KkK$, we shall write simply $\Iso(X,Y)$ and $\Iso(X)$.\par
For $X \in \KkK$ let `$\sim_{\KkK}$' be the equivalence relation on $X$ given by
$$
x \sim_{\KkK} y \iff \Phi(x) = y \textup{ for some } \Phi \in \Iso_{\KkK}(X);
$$
let $X^{(1)}$ be the quotient set $\quotient{X}{\sim_{\KkK}}$ and $\pi_X^{(1)}\dd X \to X^{(1)}$ the canonical projection.
Similarly, for any isomorphism $\Phi \in \Iso_{\KkK}(Y,Z)$ between spaces $Y, Z \in \KkK$ let $\Phi^{(1)}\dd Y^{(1)} \to
Z^{(1)}$ be the unique function such that $\pi_Z^{(1)} \circ \Phi = \Phi^{(1)} \circ \pi_Y^{(1)}$.

\begin{dfn}{disting}
Let $\KkK$ be an iso-category. By a (\textit{weak}) \textit{distinguishing} in $\KkK$
we mean any assignment $\KkK \ni X \mapsto C_X \in X^{(1)}$ such that whenever $\Phi \in \Iso_{\KkK}(K,L)$ with $K, L \in
\KkK$, then $\Phi^{(1)}(C_K) = C_L$.
\end{dfn}

More natural approach to distinguishing is the following: to each space $X \in \KkK$ assign a point $c_X \in X$
in such a way that whenever $K$ and $L$ are two isomorphic members of $\KkK$, there is an isomorphism $\Phi\dd K \to L$
which sends $c_K$ to $c_L$. However, we are interested in constructive methods of distinguishing (so, without using
the axiom of choice) and thus the original definition (\DEF{disting}) is more appropriate.\par
A very special and the most important case of distinguishing appears when the distinguished equivalence class $C_K$
consists of a single point and then we may consider $C_K$ as an element of $K$. To make this precise, we put

\begin{dfn}{str_disting}
By a \textit{strict distinguishing} in an iso-category $\KkK$ we mean any assignment $\KkK \ni X \mapsto c_X \in K$ such
that $\Phi(c_K) = c_L$ for any $K,L \in \KkK$ and each $\Phi \in \Iso_{\KkK}(K,L)$.
\end{dfn}

Strict distinguishings appear very rarely in mathematics, which the following immediate result witnesses to

\begin{pro}{cfp}
If $\KkK \ni X \mapsto c_X \in X$ is a strict distinguishing in an iso-category $\KkK$, then for every $K \in \KkK$,
$c_K$ is a common fixed point for the group $\Iso_{\KkK}(K)$. That is, $\Phi(c_K) = c_K$ for all $\Phi \in \Iso_{\KkK}(K)$.
\end{pro}

Since there are iso-categories $\KkK$ (even among those of nonempty compact spaces) in which for some spaces $K \in \KkK$
the group $\Iso_{\KkK}(K)$ has no common fixed point, a strict distinguishing is not always possible. In the next section
we introduce an iso-category for which the latter is realizable.

\SECT{Weakly convex compact metric spaces}

In the literature there are two main approaches to the notion of convexity in metric spaces. The first is related to
joining points by line segments, the second relies on generalization of the notion of the middle point between two points
by describing its global position in the space. For example, Takahashi \cite{taka} calls a metric space $(X,d)$ convex iff
for any $x, y \in X$ and every $\lambda \in (0,1)$ there is a point $z_{\lambda} \in X$ such that
\begin{equation}\label{eqn:aux1}
d(z_{\lambda},w) \leqsl (1-\lambda) d(x,w) + \lambda d(y,w)
\end{equation}
for all $w \in X$; Kijima \cite{kiji} and Yang and Zhang \cite{y-z} speak about convexity when \eqref{eqn:aux1}
with $\lambda = \frac12$ is fulfilled; while Kindler \cite{kindler} says about $\varphi$-convexity for any concave,
nondecreasing in both variables function $\varphi$ such that $\varphi(x,y) < \max(x,y)$ whenever $x \neq y$. The reader
interested in this topic is referred to the original papers of the above mentioned authors. Below we introduce
the so-called weakly convex metric spaces, the class of which includes all known to us convex metric spaces defined
by generalizing the notion of the middle point.

\begin{dfn}{weak}
A metric space $(X,d)$ is said to be \textit{weakly convex} iff for any two points $x$ and $y$ of $X$ there is a point
$z \in X$ such that for each $w \in X$:
\begin{enumerate}[(C1)]
\item $d(z,w) \leqsl \max(d(x,w),d(y,w))$,
\item $d(x,w) = d(y,w)$ provided $d(z,w) = \max(d(x,w),d(y,w))$.
\end{enumerate}
Every point $z \in X$ which satisfies (C1) and (C2) for fixed $x, y \in X$ and all $w \in X$ is said to be
a \textit{weakly middle point between $x$ and $y$}.
\end{dfn}

The reader will easily check that if $X$ is a convex subset of a normed linear space, then $X$ is weakly convex
in the sense of \DEF{weak} (for $x$ and $y \in X$ the point $z = \frac{x + y}{2}$ satisfies the conditions (C1)
and (C2)). It may also be shown that every $\RRR$-tree is weakly convex. (A complete metric space $T$ is said to be
an \textit{$\RRR$-tree} if for any two distinct points $x$ and $y$ of $T$ there is a unique homeomorphic copy
$\gamma_{x,y}$ of the interval $[0,1]$ which joins $x$ and $y$, and $\gamma_{x,y}$ is isometric to a line segment;
cf. \cite{kirk}.)\par
Our aim is to construct a strict distinguishing in the class $\WwW\CcC\CcC$ of all nonempty weakly convex compact metric
spaces (where the category is determined only by metrics). As a corollary, we shall obtain a theorem on common fixed
points in weakly convex compact metric spaces.\par
Let $(X,d)$ be a weakly convex metric space. For each $x, y \in X$ let $M(x,y)$ be the set of all weakly middle points
between $x$ and $y$ in $X$. A subset $A$ of $X$ is said to be a \textit{fully convex subspace} of $X$ iff $M(a,b) \subset
A$ for all $a, b \in A$. It is clear that a fully convex subspace of a weakly convex metric space is itself a weakly
convex metric space as well.\par
Now we shall recall the classical atributes of a metric space (see e.g. \cite{e-k} or \cite{kindler}). By $\delta(X)$
we denote the \textit{diameter} of a metric space $(X,d)$, that is, $\delta(X) := \sup_{x,y \in X} d(x,y) \in [0,+\infty]$
provided $X$ is nonempty and $\delta(\varempty) := 0$. For each $x \in X$ let $r_X(x) := \sup_{y \in X} d(x,y)$ and let
$r(X) := \inf_{x \in X} r_X(x)$ ($r(\varempty) := 0$). The number $r(X) \in [0,+\infty]$ is called the \textit{Chebyshev
radius} of $X$. Finally, the \textit{Chebyshev center} of $X$ is the set $C(X) := \{x \in X\dd\ r_X(x) = r(X)\}$.
If the latter set consists of a single point, the unique element of $C(X)$ is also called the \textit{Chebyshev center}
of $X$. The classical result states that $C(K)$ is a singleton provided $K$ is a nonempty compact convex subset
of a strictly convex normed linear space. If the assumption of strict convexity of the norm is relaxed, the set $C(K)$
may be infinite. However, $C(X)$ is nonempty for every nonempty compact metric space $X$.\par
In order to define the generalized Chebyshev center (as a uniquely determined point of a space), we introduce the following

\begin{dfn}{iterate}
The \textit{$n$-th Chebyshev center}, $C^n(X)$, of a metric space $X$ is given by the recursive formula: $C^0(X) := X$
and $C^n(X) := C(C^{n-1}(X))$ for $n > 0$. Additionally, let $C^{\infty}(X) := \bigcap_{n=0}^{\infty} C^n(X)$.
\end{dfn}

Our goal is to show that $C^{\infty}(X)$ consists of a single point provided $X$ is a nonempty weakly convex compact metric
space. To show this, we need the next result. It was proved (in a different way) in special cases by Takahashi \cite{taka}
and Kindler \cite{kindler}.

\begin{lem}{diam}
If $(X,d)$ is a weakly convex compact metric space having more than one point, then $r(X) < \delta(X)$.
\end{lem}
\begin{proof}
By an induction argument one easily proves that for every $n \geqsl 1$ and any $x_1,\ldots,x_n \in X$ there is a point
$z \in X$ such that for each $w \in X$,
\begin{enumerate}[(CC1)]
\item $d(z,w) \leqsl \max(d(x_1,w),\ldots,d(x_n,w))$,
\item $d(x_1,w) = \ldots = d(x_n,w)$ provided (CC1) is fulfilled with the equality sign.
\end{enumerate}
Now suppose, for the contrary, that $r(X) = \delta(X)$. By the compactness, there is a maximal finite system
$x_1,\ldots,x_n$ of elements of $X$ such that $d(x_j,x_k) = \delta(X)$ whenever $j \neq k$. Let $z \in X$ be a point
satisfying (CC1) and (CC2) for $x_1,\ldots,x_n$. By our assumption, $r_X(z) = \delta(X)$ and hence there is $w \in X$
such that $d(z,w) = \delta(X)$. But then, by (CC2), $d(x_1,w) = \ldots = d(x_n,w) = \delta(X)$ which denies the maximality
of the system $x_1,\ldots,x_n$.
\end{proof}

\begin{pro}{single}
For every nonempty weakly convex compact metric space $X$ the set $C^{\infty}(X)$ consists of a single point.
\end{pro}
\begin{proof}
By the compactness of $X$ and the closedness of all $C^n(X)$'s, $C^{\infty}(X)$ is nonempty and compact. It is easy
to check that $C^1(Y)$ is a fully convex subspace of $Y$ for any weakly convex metric space $Y$. This yields that $C^n(X)$
for each natural $n$, and thus $C^{\infty}(X)$ as well, is a fully convex subspace of $X$. Take $z \in C^{\infty}(X)$.
For every natural $n$, $z$ belongs to $C(C^n(X))$ and hence there is $y_n \in C^n(X)$ for which $d(z,y_n) = r(C^n(X))$.
By the definitions of the Chebyshev center and the Chebyshev radius, $r(Y) \geqsl \delta(C(Y))$ for every metric space $Y$.
We infer from this that $r(C^n(X)) \geqsl \delta(C^{n+1}(X))$ which implies that
\begin{equation}\label{eqn:aux3}
d(z,y_n) \geqsl \delta(C^{\infty}(X)) \qquad (n \in \NNN).
\end{equation}
Now let $(y_{n_k})_{k=1}^{\infty}$ be a subsequence of $(y_n)_{n=1}^{\infty}$ which converges to some $y \in X$. Then
$y \in C^{\infty}(X)$ and hence $d(z,y) = \delta(C^{\infty}(X))$, thanks to \eqref{eqn:aux3}. This shows that
$r_{C^{\infty}(X)}(z) = \delta(C^{\infty}(X))$ for every $z \in C^{\infty}(X)$ and thus $r(C^{\infty}(X)) =
\delta(C^{\infty}(X))$. Now it suffices to apply \LEM{diam} to finish the proof.
\end{proof}

\begin{dfn}{gener}
Let $X$ be a nonempty weakly convex compact metric space. The unique point of $C^{\infty}(X)$ is called
the \textit{generalized Chebyshev center} of $X$ and is denoted by $c_X$.
\end{dfn}

The construction of the generalized Chebyshev center immediately gives

\begin{thm}{strict}
Let $\WwW\CcC\CcC$ be the class of all nonempty weakly convex compact metric spaces. The assignment $\WwW\CcC\CcC \ni X
\mapsto c_X \in X$ is a strict distinguishing.
\end{thm}

The above result combined with \PRO{cfp} yields

\begin{cor}{cfp}
Let $X \in \WwW\CcC\CcC$. For every isometry $\Phi$ of $X$ onto $X$, $\Phi(c_X) = c_X$.
\end{cor}

When $X$ is a convex subset of a normed linear space, \COR{cfp} is a special case of Kakutani's fixed point theorem
on equicontinuous group of affine transformations (\cite{kak}; or \cite{rudin}). The proof presented here is constructive.
However, it works only for the specific group --- the isometry one.

\begin{rem}{affine}
The problem whether every (bijective) isometry between convex subsets of normed linear spaces is affine seems to be still
open. (Beside the classical Mazur-Ulam theorem (\cite{m-u}; or \cite[14.1]{b-l}), the author knows only one general result
\cite{mankiewicz} in this direction.) If there was a compact convex subset in a normed linear space admitting
a non-affine isometry, then \COR{cfp} would be stronger than Kakutani's fixed point theorem (in this specific case).
\end{rem}

\SECT{Central measure}

In this section we apply the results of the previous part to distinguish a measure on a compact metric space. To do this,
let us fix a nonempty compact metric space $(X,d)$. Denote by $\Prob(X)$ the set of all probabilistic Borel measures
on $X$. Equip $\Prob(X)$ with the metric $\widehat{d}$ given by the formula
\begin{equation}\label{eqn:k-r}
\widehat{d}(\mu,\nu) = \sup\Bigl\{\Bigl|\int_X f \dint{\mu} - \int_X f \dint{\nu}\Bigr|\dd\ f \in \Contr(X,\RRR)\Bigr\}
\end{equation}
where $\mu, \nu \in \Prob(X)$ and $\Contr(X,\RRR)$ stands for the family of all $d$-nonexpansive maps of $X$ into $\RRR$.
The metric $\widehat{d}$ is called the \textit{Kantorovich} (or \textit{Kantorovich-Rubenstein},
cf. \cite[Definition~2.3.1]{weaver}) metric induced by $d$. The space $(\Prob(X),\widehat{d}\,)$ is compact
and $\widehat{d}$ induces on $\Prob(X)$ the topology inherited, thanks to the Riesz characterization theorem,
from the weak-* topology of the dual Banach space of $\CCc(X,\RRR)$. It may be easily shown that $(\Prob(X),\widehat{d}\,)$
is affinely isometric to a convex subset of a normed space. Therefore $\Prob(X)$ is weakly convex. We may now introduce

\begin{dfn}{muX}
The generalized Chebyshev center of $(\Prob(X),\widehat{d}\,)$ is called the \textit{central measure} of $X$ and it is
denoted by $\mu_X$.
\end{dfn}

Now notice that every isometry $\Phi\dd X \to X$ induces an affine isometry $\widehat{\Phi}\dd \Prob(X) \to \Prob(X)$ given
by the formula $\widehat{\Phi}(\mu) = \mu \circ \Phi^{-1}$ where $\mu \circ \Phi^{-1}$ denotes the transport of the measure
$\mu \in \Prob(X)$ under the transformation $\Phi$ (that is, $(\mu \circ \Phi^{-1})(A) = \mu(\Phi^{-1}(A))$). We conclude
from \COR{cfp} that $\widehat{\Phi}(\mu_X) = \mu_X$ for all $\Phi \in \Iso(X)$; that is, $\mu_X$ is an invariant measure
for the isometry group of $X$. Again, we have obtained a constructive proof that the isometry group of an arbitrary
(nonempty) compact metric space admits an invariant measure.\par
Now suppose that $\Iso(X)$ acts transitively on $X$, i.e. for each two points $x$ and $y$ of $X$ there is $\Phi \in
\Iso(X)$ with $\Phi(x) = y$. It is known that in that case there is a unique measure invariant under every isometry of $X$
(see e.g. \cite[Theorem~2.5]{pn}). So, we get

\begin{pro}{trans}
If the isometry group of $X$ acts transitively on $X$, $\mu_X$ is the unique measure invariant under every isometry of $X$.
\end{pro}

By a \textit{metric} group we mean a metrizable topological group equipped with a left-invariant metric inducing
the topology of the group (there exists one, see e.g. \cite{berb}). As a special case of \PRO{trans} we obtain

\begin{cor}{Haar}
Let $G$ be a compact metric group. The central measure of $G$ is the Haar measure of $G$.
\end{cor}

\COR{Haar} provides a new constructive proof of the Haar measure theorem for metrizable compact groups.\par
Although in compact metric spaces with transitive actions of the isometry groups the central measures may be found thanks
of theirs very specific properties, unfortunately computing a central measure in general is very complicated. For example,
we do not know the one of $[0,1]$. The reader interested in this problem may try first to compute the central measure
of a three-point space.\\
\nl
\textbf{Question.} Is $\mu_{[0,1]}$ the Lebesgue measure?

\SECT{Quasi-nilpotent compact metric spaces}

Although we do not know whether the central measure of the unit interval is the Lebesgue measure, we are able to make
another distinguishing of measures in a special class of compact metric spaces in such a way that the distinguished measure
for the unit interval will be the Lebesgue measure. This will be done in this section.\par
Recall (see Section 2) that for a metric space $(X,d)$ the set $X^{(1)}$ is the set of all orbits of points of $X$ under
the natural action of the isometry group of $X$. It turns out that $X^{(1)}$ may be topologized by an `axiomatically'
defined metric when $(X,d)$ is compact. Precisely, we denote by $d^{(1)}$ the greatest pseudometric on $X^{(1)}$ which
makes the canonical projection $\pi_X^{(1)}\dd (X,d) \to (X^{(1)},d^{(1)})$ nonexpansive. For an arbitrary metric space
$(X,d)$, $d^{(1)}$ may not be a metric. However, we have

\begin{pro}{d1}
For every compact metric space $(X,d)$, $d^{(1)}$ is a metric on $X^{(1)}$. Moreover, for each $x, y \in X$,
\begin{multline}\label{eqn:d1}
d^{(1)}(\pi_X^{(1)}(x),\pi_X^{(1)}(y)) = \sup \{|f(\pi_X^{(1)}(x)) - f(\pi_X^{(1)}(y))|\dd\\f\dd X^{(1)} \to \RRR,\
f \circ \pi_X^{(1)} \textup{ is $d$-nonexpansive}\}.
\end{multline}
\end{pro}
\begin{proof}
The verification of \eqref{eqn:d1} is left as a simple exercise. We shall only show that $d^{(1)}$ is indeed a metric.
We shall do this with use of the variation of the Gromov-Hausdorff metric \cite{gromov,gromov2} (see also \cite{pet}).
Namely, for $a$ and $b$ in $X$ let $\varrho(a,b)$ be the least upper bound of numbers $p_H(X_1,X_2) + p((a,1),(b,2))$
where $X_j = X \times \{j\}$, $p$ is a semimetric on $X_1 \cup X_2$ such that $p((x,j),(y,j)) = d(x,y)$ for any $x, y \in
X$ and $j \in \{1,2\}$, and $p_H$ is the Hausdorff distance induced by $p$. (In other words, $\varrho(a,b)$ is
a counterpart of the Gromov-Hausdorff distance for pointed metric spaces $(X,a)$ and $(X,b)$.) As in case of the classical
Gromov-Hausdorff distance one shows that $\varrho$ is a semimetric on $X$ such that $\varrho(a,b) = 0$ iff the pointed
metric spaces $(X,a)$ and $(X,b)$ are isometric, i.e. if $\Phi(a) = b$ for some $\Phi \in \Iso(X)$. Thus $\varrho$ induces
a \textbf{metric} $\varrho^*$ on $X^{(1)}$ in such a way that $\varrho^*(\pi_X^{(1)}(x),\pi_X^{(1)}(y)) = \varrho(x,y)$
for all $x, y \in X$. Since $\varrho \leqsl d$ (because for $p((x,1),(y,2)) := d(x,y)$ one obtains $p_H(X_1,X_2) = 0$
and $p(x,y) = d(x,y)$), $\pi_X^{(1)}$ is nonexpansive with respect to the metrics $d$ and $\varrho^*$, and thus $d^{(1)}
\geqsl \varrho^*$.
\end{proof}

By \PRO{d1}, $(X,d)^{(1)} := (X^{(1)},d^{(1)})$ is a compact metric space provided $(X,d)$ is so. Thus we may repeat this
construction to obtain subsequent spaces $X^{(2)}$, $X^{(3)}$ and so on. Namely, for a compact metric space let
$(X^{(0)},d^{(0)}) = (X,d)$ and $(X^{(n)},d^{(n)}) = (X^{(n-1)},d^{(n-1)})^{(1)}$ for $n > 0$. Notice that $\delta(X^{(n)})
\leqsl \delta(X^{(n-1)})$ and thus the sequence $(\delta(X^{(n)})_{n=1}^{\infty}$ is convergent. We introduce the following

\begin{dfn}{quasi}
A compact metric space $X$ is said to be \textit{quasi-nilpotent} iff $\lim_{n\to\infty} \delta(X^{(n)}) = 0$.
\end{dfn}

The class of quasi-nilpotent compact metric spaces includes all spaces on which their isometry groups acts transitively.
One may think that such spaces have to have rich isometry groups. The next example shows that it is not the rule.

\begin{exm}{01}
Let $(X,d)$ be the interval $[a,b]$ with the natural metric. Observe that the isometry group of $X$ is very poor ---
there is only one isometry on $X$ different from the identity map. However, $X$ is quasi-nilpotent. To see this,
it suffices to show that $X^{(1)}$ is isometric to $[a/2,b/2]$. But this may easily be shown by means of \eqref{eqn:d1}.
\end{exm}

Now we shall distinguish a special measure on a quasi-nilpotent (nonempty) compact metric space, which may also be called
central. For a convex subset $K$ of a normed linear space let $\Fix(K)$ be the set of all fixed points under every affine
isometry of $K$ onto $K$. The set $K$ is convex as well (however, it may be empty). Further, let $\Fix^0(K) := K$ and for
natural $n > 0$ let $\Fix^n(K) = \Fix(\Fix^{n-1}(K))$. Finally, put $\Fix^{\infty}(K) = \bigcap_{n=0}^{\infty} \Fix^n(K)$.
Note that $\Fix^{\infty}(K)$ is convex and if $K$ is compact, $\Fix^{\infty}(K)$ is nonempty.\par
Now let $X$ be a nonempty compact metric space and $\Prob(X)$ be equipped with the Kantorovich metric induced by the metric
of $X$. Let $\Delta(X) = \Fix^{\infty}(\Prob(X))$. In the sequel we shall prove that $\Delta(X)$ consists of a single
measure iff $X$ is quasi-nilpotent. In fact, this follows from the following

\begin{thm}{iso}
For a nonempty compact metric space $(X,d)$ the function
$$
\Psi\dd \Fix(\Prob(X)) \ni \mu \mapsto \mu \circ \pi_X^{-1} \in \Prob(X^{(1)})
$$
is an affine isometry of $\Fix(\Prob(X))$ onto $\Prob(X^{(1)})$. In particular, $\delta(\Fix(\Prob(X)) = \delta(X^{(1)})$.
\end{thm}
\begin{proof}
Since $\delta(\Prob(Y)) = \delta(Y)$ for every compact metric space $Y$, it suffices to prove the first assertion.
By \cite{pn}, $\Psi$ is an affine bijection. So, we only need to check that $\Psi$ is isometric. Fix $\mu_1, \mu_2 \in
\Fix(\Prob(X))$ and put $\nu_j = \Psi(\mu_j)$. If $u \in \Contr(X^{(1)},\RRR)$, then $\int_{X^{(1)}} u \dint{\nu_j}
= \int_X u \circ \pi_X^{(1)} \dint{\mu_j}$. This, combined with \eqref{eqn:k-r}, gives $\widehat{d}(\mu_1,\mu_2) \geqsl
\widehat{d^{(1)}}(\nu_1,\nu_2)$. Conversely, if $v \in \Contr(X,\RRR)$, then, since $\mu_j \in \Fix(\Prob(X))$,
$\int_X v \dint{\mu_j} = \int_X v \circ \Phi \dint{\mu_j}$ for every $\Phi \in \Iso(X)$.
Now by \cite[Proposition~2.5]{pn2}, the closed convex hull (in the topology of uniform convergence) of the set
$\{v \circ \Phi\dd\ \Phi \in \Iso(X)\}$ contains a map $w\dd X \to \RRR$ such that
\begin{equation}\label{eqn:aux5}
w \circ \Phi = w
\end{equation}
for all $\Phi \in \Iso(X)$. This implies that $w \in \Contr(X,\RRR)$ and $\int_X v \dint{\mu_j} = \int_X w \dint{\mu_j}$.
We infer from \eqref{eqn:aux5} that there is $w_0\dd X^{(1)} \to \RRR$ such that $w = w_0 \circ \pi_X^{(1)}$. The latter
connection and \eqref{eqn:d1} yield that $w_0 \in \Contr(X^{(1)},\RRR)$. So, we finally obtain $|\int_X v \dint{\mu_1}
- \int_X v \dint{\mu_2}| = |\int_X w_0 \circ \pi_X^{(1)} \dint{\mu_1} - \int_X w_0 \circ \pi_X^{(1)} \dint{\mu_2}|
= |\int_{X^{(1)}} w_0 \dint{\nu_1} - \int_{X^{(1)}} w_0 \dint{\nu_2}| \leqsl \widehat{d^{(1)}}(\nu_1,\nu_2)$,
which finishes the proof.
\end{proof}

Now \THM{iso} and induction argument give

\begin{pro}{dn}
If $X$ is a nonempty compact metric space, then $\delta(\Fix^n(\Prob(X))) = \delta(X^{(n)})$ for each natural $n$.
\end{pro}

\begin{cor}{quasi}
Let $X$ be a nonempty compact metric space. $\Delta(X)$ consists of a single measure iff $X$ is quasi-nilpotent.
\end{cor}

\begin{dfn}{lambda}
Let $X$ be a nonempty quasi-nilpotent compact metric space. The unique member of $\Delta(X)$ is denoted by $\lambda_X$
and it is called the \textit{central measure of $X$ of a second kind}.
\end{dfn}

Since $\lambda_X \in \Fix(\Prob(X))$, the central measure of $X$ of a second kind is invariant under every isometry
of $X$. We conclude from this that $\lambda_X = \mu_X$ provided $X$ is a compact metric space such that $X^{(1)}$ is
a singleton (i.e. if the isometry group of $X$ acts transitively on $X$). We end the section with

\begin{pro}{Lebesgue}
The central measure of $[0,1]$ of a second kind coincides with the Lebesgue measure on $[0,1]$.
\end{pro}
\begin{proof}
For simplicity, put $I = [0,1]$ and $\lambda = \lambda_I$. We infer from the relation $\lambda \in \Fix(\Prob(I))$ that
$\lambda = \lambda \circ u_1$ where $u_1\dd I \ni t \mapsto |t - 1/2| \in [0,1/2]$. Similarly, since $\lambda \in
\Fix^n(\Prob(I))$, $\lambda$ is invariant under the map $u_n\dd [0,1/2^{n-1}] \ni t \mapsto |t - 1/2^n| \in [0,1/2^n]$.
One deduces from this that $\lambda(\{\frac{k}{2^n}\}) = 0$ and $\lambda([\frac{k-1}{2^n},\frac{k}{2^n}]) = \frac{1}{2^n}$
for each natural $k$ and $n$ with $1 \leqsl k \leqsl 2^n$, and hence $\lambda$ is the Lebesgue measure. The details are
left for the reader.
\end{proof}

\SECT{Distinguishing dense subsets}

We know that strict distinguishing is impossible in general. However, one may still ask whether it is possible to define
a distinguishing in the class $\KkK$ of all nonempty compact metric spaces. This part is devoted to the solution of this
problem. We shall show that there is a sequence of distinguishings $\KkK \ni K \mapsto C_n(K) \in K^{(1)}\ (n \geqsl 1)$
such that for every $K \in \KkK$, whenever $c_n \in K$ satisfies $\pi_K^{(1)}(c_n) = C_n(K)$, then the set $\{c_n\dd\
n \geqsl 1\}$ is dense in $K$.\par
Let $\NNN = \{0,1,2,\ldots\}$. Fix an infinite compact metric space $(X,d)$. Instead of constructing an `intrinsic' dense
subset of $X$, we shall construct a metric $\varrho_X$ on $\NNN$ such that $(\NNN,\varrho_X)$ is isometric to a dense
subset of $X$. Suppose for some $n \in \NNN$ we have defined the metric $\varrho_X$ on $\{0,\ldots,n\}$ in such a way that
the space $(\{0,\ldots,n\},\varrho_X)$ is isometrically embeddable into $X$ (for $n = 0$ we have nothing to do). Put
\begin{equation}\label{eqn:Fn}
F_n := \{(x_0,\ldots,x_n) \in X^{n+1}\dd\ d(x_j,x_k) = \varrho_X(j,k),\ j,k=0,\ldots,n\}
\end{equation}
and
\begin{equation}\label{eqn:fn}
f_n\dd F_n \times X \ni (x_0,\ldots,x_n;x) \mapsto \min(d(x_0,x),\ldots,d(x_n,x)) \in \RRR.
\end{equation}
By our assumption, $F_n$ is nonempty. Next, let
\begin{equation}\label{eqn:A0n}
A_0^{n+1} := \{(x;y) \in F_n \times X\dd\ f_n(x;y) = \max f_n(F_n \times X)\}.
\end{equation}
Now inductively define sets $A_j^{n+1}$ for $j=1,\ldots,n+1$ by
\begin{multline}\label{eqn:Ajn}
A_j^{n+1} := \Bigl\{(y_0,\ldots,y_n;y) \in A_{j-1}^{n+1}\dd\\ d(y_{j-1},y) = \max \{d(x_{j-1},x)\bigr|\
(x_0,\ldots,x_n;x) \in A_{j-1}^{n+1}\}\Bigr\}.
\end{multline}
The reader will easily check (by induction and the compactness argument) that each of the sets $A_j^{n+1}$'s is nonempty.
Now take an arbitrary $(x_0,\ldots,x_n;x_{n+1}) \in A_{n+1}^{n+1}$ and put $\varrho_X(j,n+1) := d(x_j,x_{n+1})$
for $j=0,\ldots,n+1$. Observe that this definition is independent of the choice of $(x_0,\ldots,x_n;x) \in A_{n+1}^{n+1}$.
It is also clear that $\varrho_X$ is a metric (not only a semimetric) on $\{0,1,\ldots,n+1\}$ (because $X$ is
infinite).\par
In this way we obtain a metric $\varrho_X$ on $\NNN$ such that $\varrho_X = \varrho_Y$ for every space $Y$ isometric
to $X$. We claim that

\begin{pro}{dense}
For every infinite compact metric space $(X,d)$, the space $(\NNN,\varrho_X)$ is isometric to a dense subset of $X$.
\end{pro}
\begin{proof}
For each $n \in \NNN$ let $P_n$ be the set of all sequences $(x_m)_{m=0}^{\infty} \in X^{\NNN}$ such that the function
$(\{0,\ldots,n\},\varrho_X) \ni j \mapsto x_j \in (X,d)$ is isometric. By construction of $\varrho_X$, $P_n$ is nonempty.
It is also clear that $P_n$ is closed in $X^{\NNN}$ and that $P_n \supset P_{n+1}$. Therefore, by the compactness of $X$,
the intersection $\bigcap_{n=0}^{\infty} P_n$ is nonempty. We infer from this that there is an isometric function
$\Phi$ of $(\NNN,\varrho_X)$ into $(X,d)$. We claim that $\Phi(\NNN)$ is dense in $X$. Suppose, for the contrary,
that there is $x \in X$ and $r > 0$ such that
\begin{equation}\label{eqn:fnx}
d(x,\Phi(n)) \geqsl r
\end{equation}
for every $n \in \NNN$. Note that $(\Phi(0),\ldots,\Phi(n+1)) \in A_{n+1}^{n+1} \subset A_0^{n+1} \subset F_n \times X$
for any $n \in \NNN$, where $A_j^{n+1}$'s and $F_n$'s are given by \eqref{eqn:Ajn}, \eqref{eqn:A0n} and \eqref{eqn:Fn}.
So, \eqref{eqn:fnx} yields $\max f_n(F_n \times X) \geqsl r$ for $f_n$'s given by \eqref{eqn:fn}. Finally, we conclude
from the relation $(\Phi(0),\ldots,\Phi(n+1)) \in A_0^{n+1}$ and \eqref{eqn:fnx} that
$f_n(\Phi(0),\ldots,\Phi(n);\Phi(n+1)) \geqsl r$ which means that $d(\Phi(j),\Phi(k)) \geqsl r$ for $j < k$.
But this denies the compactness of $X$.
\end{proof}

By a \textit{representation} of the metric $\varrho_X$ we mean any isometric function of $(\NNN,\varrho_X)$
into $(X,d)$, provided $X$ is infinite.\par
If $X$ is finite and has $n$ elements, we may repeat the above construction to obtain a metric $\varrho_X$
on $\{0,\ldots,n-1\}$ which makes this set isometric to $X$. In that case by a \textit{representation} of $\varrho_X$
we mean any function $\Phi\dd \NNN \to X$ such that $\Phi$ is isometric on $\{0,\ldots,n-1\}$ (with respect to the metrics
$\varrho_X$ and $d$) and $\Phi(k) = \Phi(n-1)$ for $k < n-1$.\par
We may ask how many representations has the metric $\varrho_X$ for an arbitrary space $X$. The answer to this gives

\begin{pro}{bij}
Let $X$ be a nonempty compact metric space and $\Phi_0\dd \NNN \to X$ a representation of $\varrho_X$. The function
$\Psi \mapsto \Psi \circ \Phi_0$ establishes a one-to-one correspondence between isometries [$\Psi$] of $X$
and representations [$\Psi \circ \Phi_0$] of $\varrho_X$.
\end{pro}

The above result is an immediate consequence of \PRO{dense} (and the well known fact that every isometric map of a compact
metric space into itself is onto \cite{lin}). \PRO{bij} says that if the isometry group of a space $X$ is poor, there are
only few representations of $\varrho_X$. In the opposite, if there are many representations, the isometry group of $X$
is rich. Both the situations are interesting.\par
Now we pass to distinguishing of points. Observe that whatever representation $\Phi\dd \NNN \to X$ of $\varrho_X$ we take,
the function $\pi_X^{(1)} \circ \Phi$ is the same. It follows from the latter that the definition
$C_n(X) := \pi_X^{(1)}(\Phi(n))$ where $n \in \NNN$ and $\Phi$ is any representation of $\varrho_X$ is correct.
We now clearly have

\begin{pro}{disting}
For each $n \in \NNN$, the assignment $\KkK \ni K \mapsto C_n(K) \in K^{(1)}$ is a distinguishing.
\end{pro}

In studying the class of separable complete metric spaces, especially in theory of random metric spaces
(cf. \cite{vershik1,vershik2}), one of methods is to consider the set of all metrics $\Dd$ on $\NNN$ and to make
the assignment $\Dd \ni d \mapsto$ `the completion of $(\NNN,d)$'. In other words, the `world' of infinite separable
complete metric spaces may be identified (by this assignment) with the `world' of metrics on $\NNN$. This is quite natural
approach, however, there is no one-to-one correspondence between the members of these two worlds. The distinguishing
of dense subsets of compact metric spaces constructed in this section may be seen as an example of the `inverse function'
to the above assignment after restricting the considerations to totally bounded metrics on $\NNN$.

\end{document}